\documentclass[review]{elsarticle}
\usepackage{setspace}
\usepackage{longtable}
\usepackage{lineno,hyperref}
\usepackage{amssymb,amsfonts,amsthm}
\usepackage{amsmath}
\usepackage{graphicx}
\usepackage{float}
\usepackage{caption}
\usepackage{subcaption}
\usepackage{parselines}
\usepackage{bm}
\usepackage{multirow}
\pdfoutput=1

\usepackage[margin=1.0in]{geometry}
\begin{document}

\title{Neural Network with Local Converging Input (NNLCI) for Supersonic Flow Problems with Unstructured Grids}
\author{Weiming Ding\footnotemark[1], Haoxiang Huang\footnotemark[2],\; Tzu Jung Lee\footnotemark[3],\; Yingjie Liu\footnotemark[4] \; and Vigor Yang\footnotemark[5]}

\begin{abstract}
In recent years, surrogate models based on deep neural networks (DNN) have been widely used to solve partial differential equations, which were traditionally handled by means of numerical simulations.  This kind of surrogate models, however, focuses on global interpolation of the training dataset, and thus requires a large network structure.  The process is both time consuming and computationally costly, thereby restricting their use for high-fidelity prediction of complex physical problems. In the present study, we develop a neural network with local converging input (NNLCI) for high-fidelity prediction using unstructured data. The framework utilizes the local domain of dependence with converging coarse solutions as input, which greatly reduces computational resource and training time. As a validation case, the NNLCI method is applied to study inviscid supersonic flows in channels with bumps. Different bump geometries and locations are considered to benchmark  the effectiveness and versability of the proposed approach.  Detailed flow structures, including shock-wave interactions, are examined systematically.
\end{abstract}

\maketitle

\renewcommand{\thefootnote}{\fnsymbol{footnote}}
\footnotetext{Key words: neural network, neural networks with local converging inputs, physics informed machine learning, conservation laws, differential equation, multi-fidelity optimization.}
\footnotetext{1. {\tt E-mail: wmding@gatech.edu}. Daniel Guggenheim School of Aerospace Engineering, Georgia Institute of Technology, Atlanta, GA 30332, USA.}
\footnotetext{2. {\tt E-mail: hcwong@gatech.edu}. Woodruff School of Mechanical Engineering, Georgia Institute of Technology, Atlanta, GA 30332, USA.}
\footnotetext{3. {\tt E-mail: tlee452@gatech.edu}. Daniel Guggenheim School of Aerospace Engineering, Georgia Institute of Technology, Atlanta, GA 30332, USA.}
\footnotetext{4. {\tt E-mail: yingjie@math.gatech.edu}. School of Mathematics, Georgia Institute of Technology, Atlanta, GA 30332, USA.}
\footnotetext{5. {\tt E-mail: vigor.yang@aerospace.gatech.edu}. Daniel Guggenheim School of Aerospace Engineering, Georgia Institute of Technology, Atlanta, GA 30332, USA.}
\renewcommand{\thefootnote}{\arabic{footnote}}

\makeatletter
\def\ps@pprintTitle{%
  \let\@oddhead\@empty
  \let\@evenhead\@empty
  \let\@oddfoot\@empty
  \let\@evenfoot\@oddfoot
}
\makeatother

\section{Introduction}
Numerical solution of partial differential equations is an essential aspect of learning about physical phenomena in many natural and engineering science disciplines. The spatiotemporal discretization of the underlying governing equations is computationally expensive and time-consuming. In recent years, with the rapid development of machine learning techniques, researchers have proposed alternative ways to handle the problems more efficiently. Data-driven surrogate models have been developed to improve, or even replace, traditional numerical simulations by mapping between the problem setting and its solution. 
Yang and colleagues \cite{1,2,3,4,5} developed a POD-based surrogate model for emulating detailed spatio-temporally evolving flows in swirl injectors.  Jose et al. \cite{6} utilized U-Net for the parametric study of shear coaxial injector flow evolution. Milan et al. \cite{7} applied a deep neural network (DNN) to study automotive fuel injector design. On the other hand, the deep Galerkin method, deep Ritz method and physics-informed neural network (PINN) have been extensively studied to solve partial differential equations (PDE) and physical problems by enforcing the boundary and initial conditions, and the functional forms of underlying governing equations \cite{8,9,10}. The usage of neural operators further improves the ability to address various problems \cite{11,12,13}.  Researchers also implemented multi-fidelity approximation of high-resolution solutions (also known as super-resolution or upscaling), due to the limitation in computation resources or high-fidelity data. For instance, Kennedy \cite{14} proposed a co-Kriging model, which uses both the high- and low-fidelity data. It was then applied to the uncertainty quantification of beam frequency and the prediction of the airfoil pressure field\cite{15,16}. Erichson et al. \cite{17} employed shallow neural network technique to reconstruct high-resolution solutions from limited observations.

Although the above methods have shown promising results, they have limits that prohibit their application to the high-resolution prediction of complex nonlinear physical problems. Most existing DNN-based surrogate models are based on global computation of the physical field. The mapping between the low-fidelity latent space and the high-fidelity data over the entire computational domain needs to be developed, which requires a large covariance matrix or neural network structure, and is computationally intensive. In addition, in physics-based methods \cite{8,9,10}, the computation of the functional forms adds the burden on the training process of the network. This limits its application to complicated problems involving discontinuities, particularly to interactions among such discontinuities.

To address this issue, local based surrogate models were proposed. Trask et al. \cite{18} developed convolutional neural network (CNN) with a generalized moving least squares (GMLS) method. The scattered data inputs are used to construct the local regression functions. However, the local feature and underlying information of the data is not fully utilized. Huang et al. \cite{19,20} established a novel method, known as Neural Networks with Local Converging Inputs (NNLCI), to solve conservation laws at low cost. This method predicts the high-resolution solution at a space-time location from two converging, low-fidelity local input solution patches. With the use of local domain of dependence, the method extracts important local features for accurate predictions, while at the same time reducing the need of computational resources and training data. The NNLCI method has shown great prediction accuracy for solving 1D \cite{19} and 2D \cite{20} Euler equations and Maxwell’s equations \cite{21}.

The application of NNLCI to structured data is relatively easy to implement, due to its nature of local domain of dependence. Many scientific and engineering problems, however, involve irregularly structured data sets. Extension of NNLCI to such situations is necessary. In the present work, we will develop a new NNLCI method for unstructured data. As a validation case, inviscid flow through a converging-diverging channel with two smooth Gaussian bumps was studied systematically. The new model is capable of capturing the flow behaviors in the entire field, including regions with smooth variations and shock discontinuities. 

This paper is structured as follows. Section~\ref{Sec: Numerical Framework} describes the theoretical formulation and numerical scheme for inviscid channel flows. The data generation and pre-processing steps are also discussed. Section ~\ref{Sec: NNLCI} introduces the neural network with local converging inputs (NNLCI). The method for determining the local domain of dependence is developed for unstructured grids. Section ~\ref{Sec: Results and Discussion} presents the results of the proposed method.  The effectiveness of the new approach is demonstrated by a variety of channel geometries. In Section \ref{Sec: conclusion}, the conclusions are reported.

\section{Theoretical and Numerical Framework}
\label{Sec: Numerical Framework}
\subsection{Problem Setup}
\begin{figure}[ht!]
    \centering
    \includegraphics[width=4.5in]{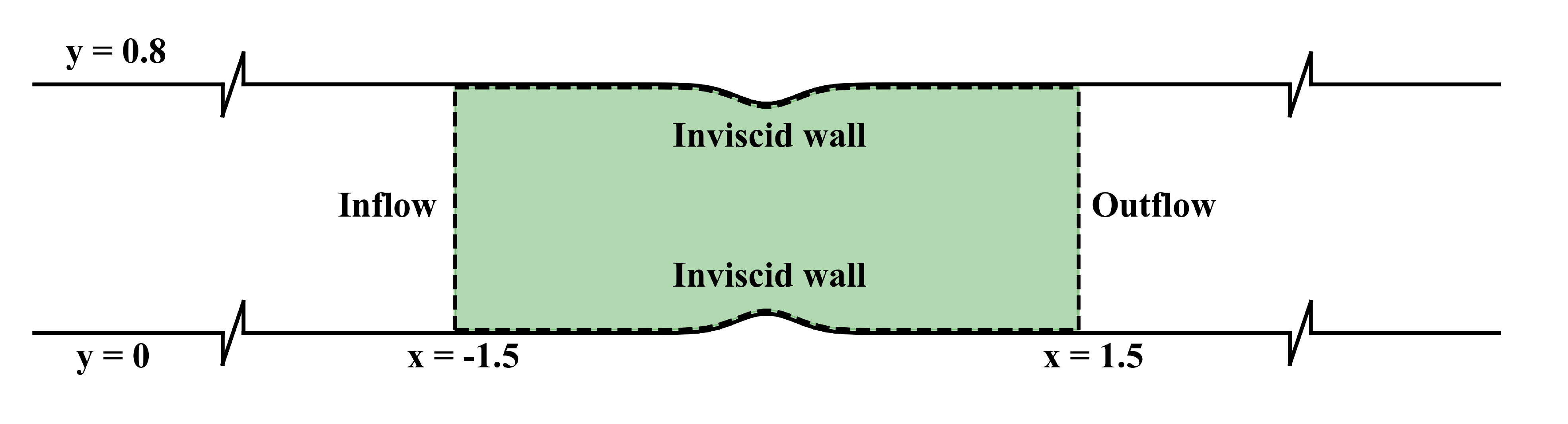}
    \label{t=0}
    \caption{Computational domain for inviscid flow through a channel with two smooth Gaussian bumps}
    \label{fig:domain}
\end{figure}

In this study, we consider a 2D inviscid flow through a channel with two smooth Gaussian bumps, as shown schematically in Fig.~\ref{fig:domain}. The computational domain is bounded by  $x \in [-1.5,1.5]$ horizontally, and $y \in [0.0,0.8]$ in the y-direction. Two Gaussian bumps are placed on the top and bottom walls of the channel. The lower bump geometry is fixed and defined by:
\begin{equation}
    y = 0.0625e^{-25x^2}.
\end{equation}
The upper bump is perturbed from the original location and is defined by:
\begin{equation}
    y = 0.8-0.0625e^{-25(x-\Delta x)^2}.
\end{equation}
where $\Delta x$ stands for the perturbation of the bump location.

The problem is governed by the 2D Euler equations for compressible flows:
\begin{equation}
    \frac{\partial \rho}{\partial t} + \frac{\partial (\rho u)}{\partial x} + \frac{\partial (\rho v)}{\partial y} = 0,
\end{equation}
\begin{equation}
    \frac{\partial (\rho u)}{\partial t} + \frac{\partial (\rho u^2+p)}{\partial x} + \frac{\partial (\rho uv)}{\partial y} = 0,
\end{equation}
\begin{equation}
    \frac{\partial (\rho v)}{\partial t} + \frac{\partial (\rho vu)}{\partial x} + \frac{\partial (\rho v^2+p)}{\partial y} = 0,
\end{equation}
\begin{equation}
    \frac{\partial (\rho E)}{\partial t} + \frac{\partial (\rho uH)}{\partial x} + \frac{\partial (\rho vH)}{\partial y} = 0,
\end{equation}
where $\rho$,$u$,$v$, $p$and $E$ are density, x-velocity, y-velocity, pressure, and total energy, respectively. $H=E+p/\rho$ is the total enthalpy. Calorically perfect gas is assumed. The resulting equation of state takes the form:
\begin{equation}
    p = (\gamma-1)(\rho E-\frac{1}{2}\rho \| \textbf{v} \|_2^2).
\end{equation}
The ratio of specific heats $\gamma$ is taken as 1.4. At the top and bottom wall, we apply the inviscid wall boundary condition. At the inflow, we specify the total temperature $T_t$ and total pressure $p_t$ as:
\begin{equation}
    \frac{T_{t}}{T_{\infty}} = 1+\frac{\gamma-1}{2}M_{\infty}^2,
\end{equation}
\begin{equation}
    \frac{p_t}{p_{\infty}} = (\frac{T_t}{T_{\infty}})^{\gamma/(\gamma-1)}.
\end{equation}
Both the inlet total temperature and total pressure depend on the inlet Mach number. For this study, we choose the inlet Mach number to be  $M_{\infty}=2.0$.

A finite-volume solver is implemented for the problem using the MUSCL scheme~\cite{22} with Rusanov flux~\cite{23}. The improved Euler, which is a total variation diminishing (TVD) second-order Runge-Kutta (RK2) scheme~\cite{24} is used for the explicit time marching. In the end, a steady state solution can be obtained. 

\subsection{Data Generation and Pre-processing}
To generate the training and testing data set, different bump location perturbations need to be considered. We translate the upper bump location in the x-direction with $\Delta x=0.00$, $\pm0.15$, $\pm0.30$, $\pm0.45$, and $ \pm0.60$. Similarly, for the testing dataset, we randomly perturb the upper bump location within the range of $[-0.60, 0.60]$. In addition, for the training cases, the free-stream Mach number is perturbed by $\pm5\%$. This ensures that the training dataset contains all the features of highly nonlinear flow behavior and increases the robustness of the network.  Table~\ref{tab:doe} shows the details of the training and testing dataset.
\begin{table}[ht!]
\renewcommand{\arraystretch}{1.5} 
\centering
\begin{tabular}{ccc}
\hline
Set & Description & $M_{\infty}$ Perturbation\\
\hline
Training & $\Delta x=0.00$, $\pm0.15$, $\pm0.30$, $\pm0.45$, and $ \pm0.60$ & $\Delta M_{\infty}=0$, $\pm5\%$\\
Testing & $\Delta x=0.12$, $-0.19$, $-0.35$, $0.44$ & None\\
\hline
\end{tabular}
\caption{\label{tab:doe} Summary of the training and testing cases. For each bump translation case, the free-stream Mach number $M_{\infty}$ is perturbed by $\pm5 \%$}
\end{table}

\begin{figure}[ht!]
    \centering
    \includegraphics[width=6in]{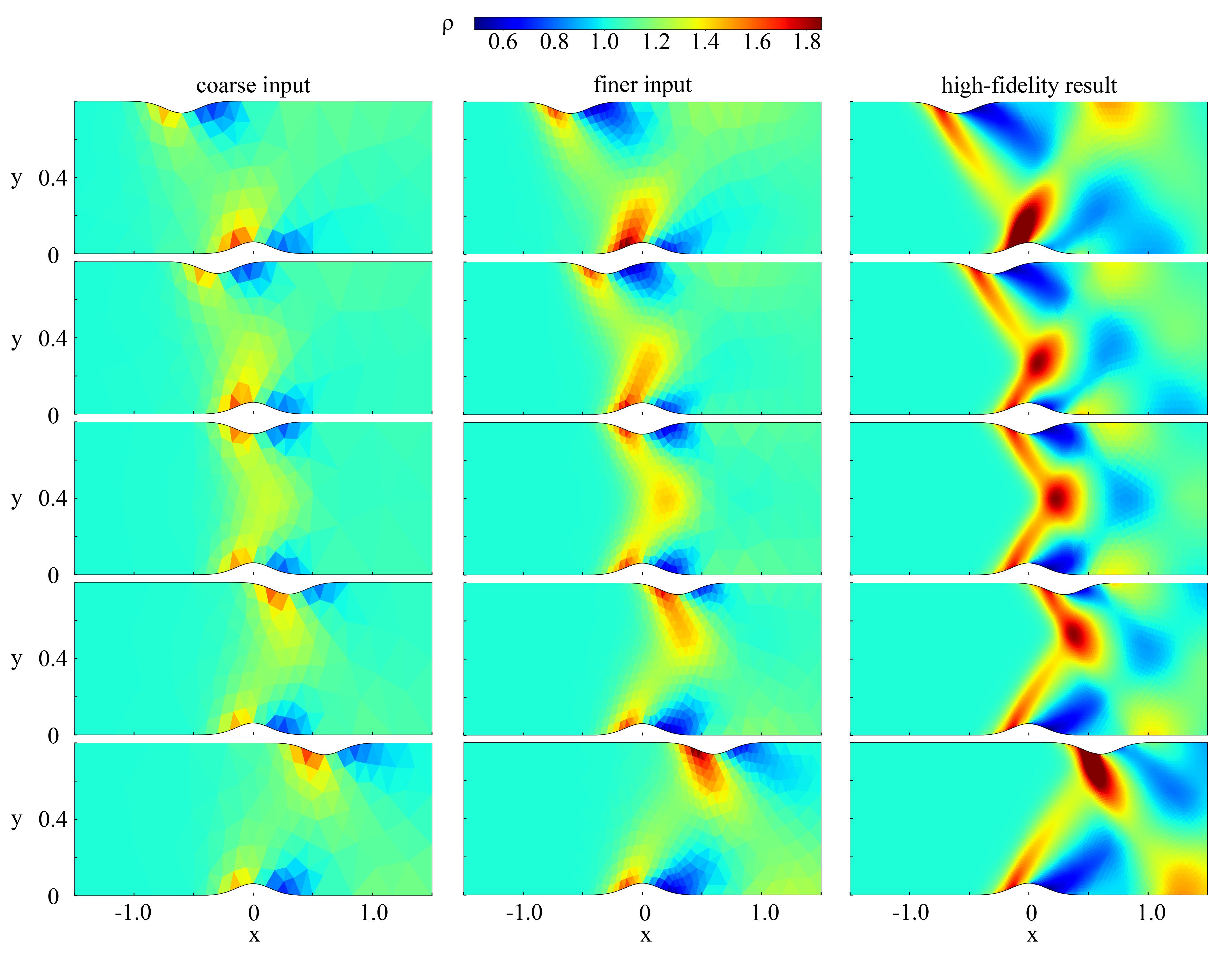}
    \caption{Example of density fields for case $\Delta x = -0.60$, $-0.30$, $0.00$, $0.30$ and $0.60$, $M_{\infty}=2.0$. The coarse, finer and high-resolution solutions are shown, from left to right. }
    \label{fig:Simulation_density}
\end{figure}

\begin{figure}[ht!]
    \centering
    \includegraphics[width=6in]{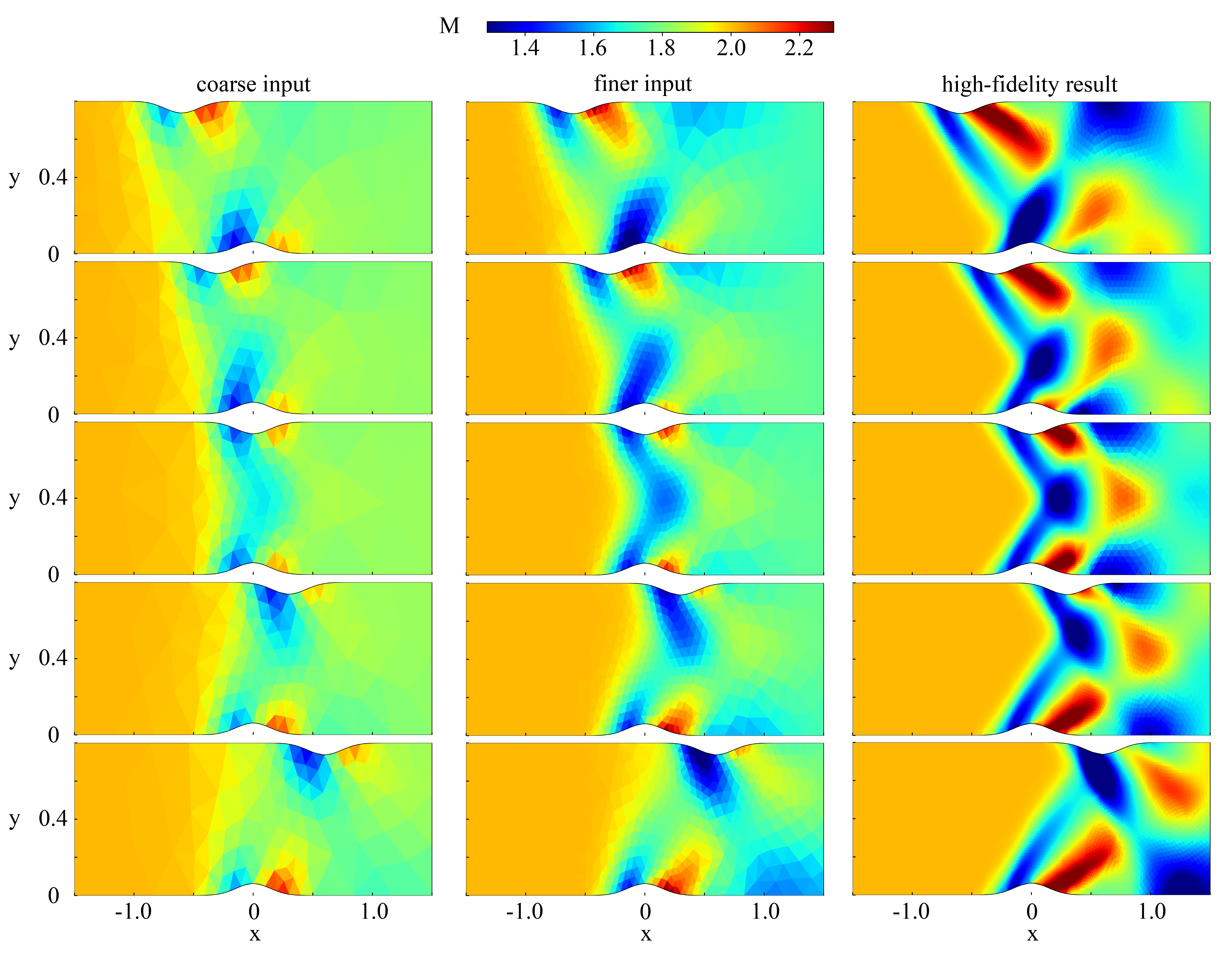}
    \caption{Calculated Mach-number fields for cases $\Delta x = -0.60$, $-0.30$, $0.00$, $0.30$ and $0.60$, $M_{\infty}=2.0$. The coarse, finer inputs, and high-fidelity results are shown, from left to right.}
    \label{fig:Simulation_Mach}
\end{figure}

For each geometry setting, we solve the problems on unstructured triangular grids with three different resolutions: coarse, finer, and high-resolution grids. To obtain the finer and high-resolution meshes, we implement uniform mesh adaptation from the coarse grid.; the numbers of cells are 200, 800, and 12800, respectively. Fig.~\ref{fig:Simulation_density} and ~\ref{fig:Simulation_Mach} show the calculated density-gradient and Mach-number fields, respectively, for the upper bump translation at  $\Delta x = -0.60$, $-0.30$, $0.00$, $0.30$ and $0.60$. For each case, the coarse and finer inputs, and high-fidelity results are shown from left to right. From observation, as the upper bump is translated, the flow behavior in the shock and downstream regions changes considerably. In addition, the high-fidelity simulation results exhibit much richer details in the shock intersection and expansion regions. 

\begin{figure}[ht!]
    \centering
    \includegraphics[width=6in]{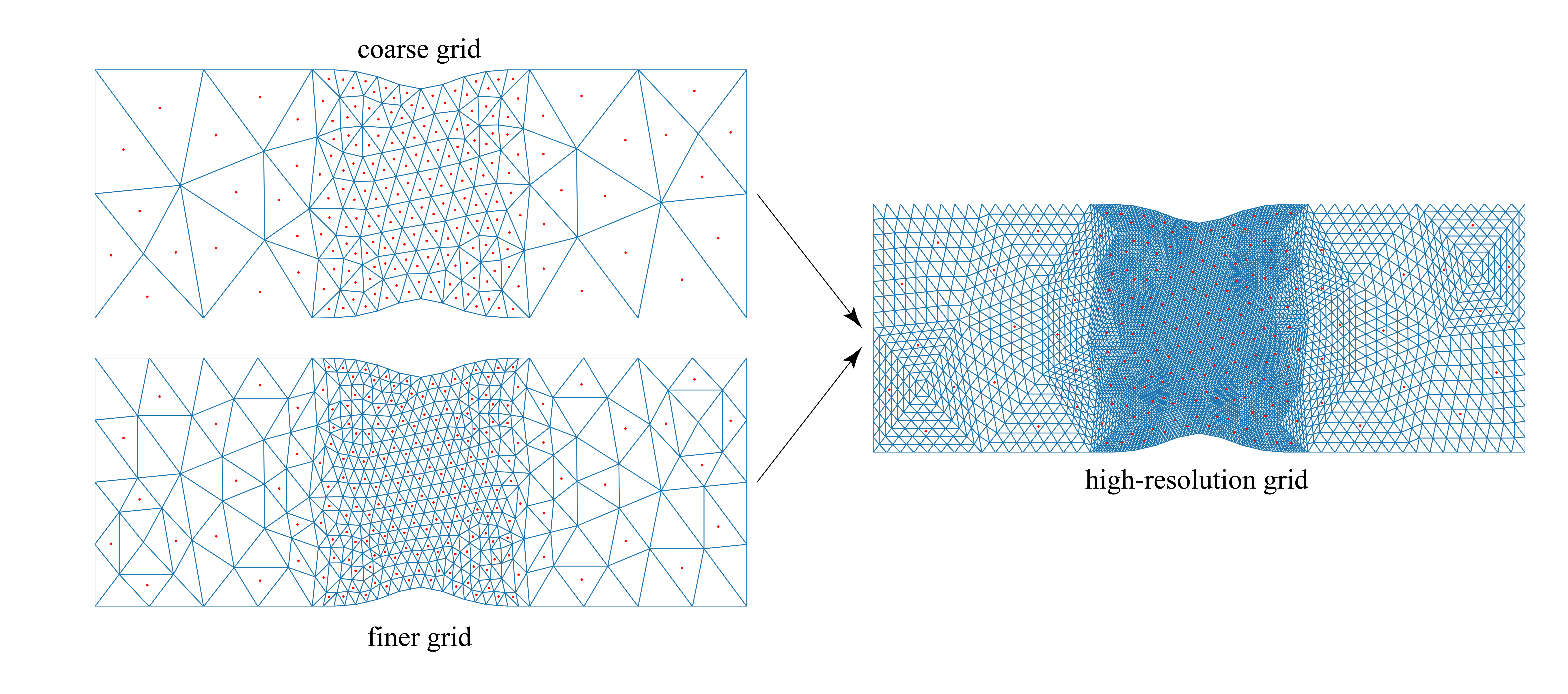}
    \caption{Example of the unstructured NNLCI training dataset. The coarse and finer data are used as the inputs, and the finest data is used as the image for training.}
    \label{fig:NNLCI_Inputs}
\end{figure}

To fully utilize the low-fidelity information and extract common flow features from different bump translation cases, special pre-processing needs to be performed on the coarse and finer grid data. Fig.~\ref{fig:NNLCI_Inputs} shows an example of the data selection process. Only the cell centroids of the coarse grids (red dots in plot) are selected as the training locations. The data at the corresponding locations in the finer and high-resolution grids are extracted to form the input-image pairs for the training dataset. That is, 200 points are used in each training case. This avoids the computationally expensive interpolation of a large dataset, and accelerates the training process.

For each data point, all four state variables $\bm{u} = [\rho, \rho u, \rho v, \rho E]$ are used. Normalization and standardization are needed to facilitate the training of the neural network. In this study, the data are rescaled by the difference between maximum and minimum values:
\begin{equation}
    \Tilde{\bm{u}} = \frac{\bm{u}-min(\bm{u})}{max(\bm{u})-min(\bm{u})}
\end{equation}
Correspondingly, the Tanh function is used as the activation function in the neural network:
\begin{equation}
    f(x) = \frac{e^x-e^{-x}}{e^x+e^{-x}}
\end{equation}

\section{Neural Network with Local Converging Inputs}
\label{Sec: NNLCI}
\begin{figure}[ht!]
    \centering
    \includegraphics[width=6in]{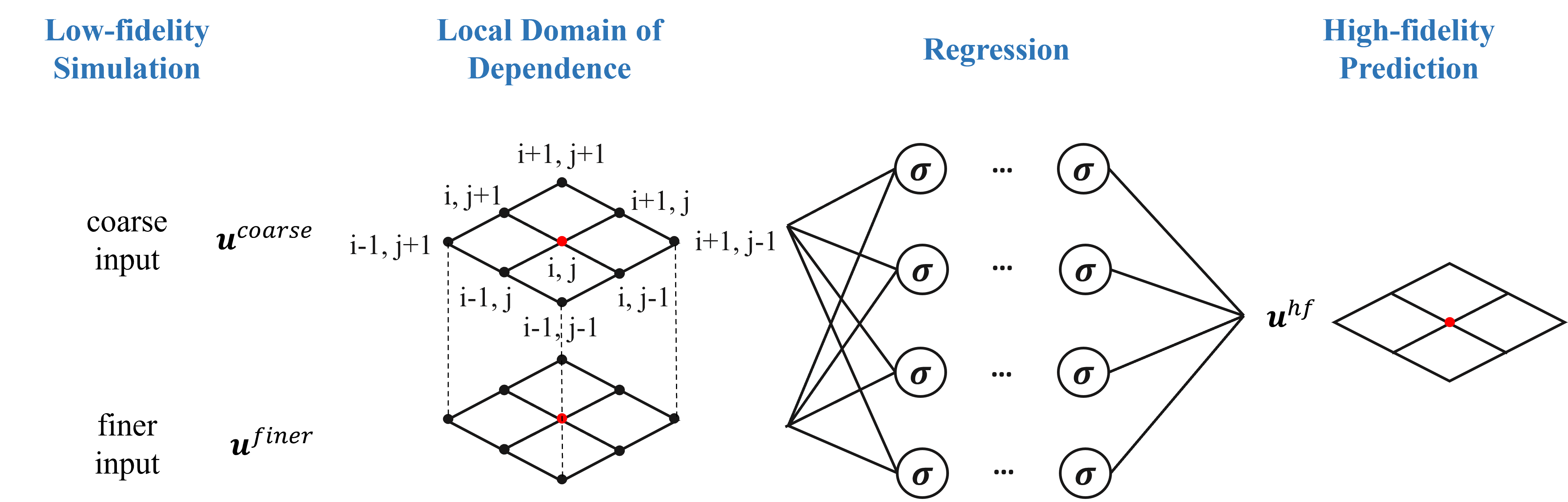}
    \label{t=0}
    \caption{Overview of the NNLCI structure. Four major steps are involved: multi-fidelity simulation, local domain of dependence, neural network regression and high-fidelity prediction.}
    \label{fig:Network_Structure}
\end{figure}

In this section, we introduce the setup of Neural Network with Local Converging Inputs (NNLCI) for unstructured data. Fig.~\ref{fig:Network_Structure} provides an overview of the NNLCI.

\begin{figure}[ht!]
    \centering
    \includegraphics[width=6in]{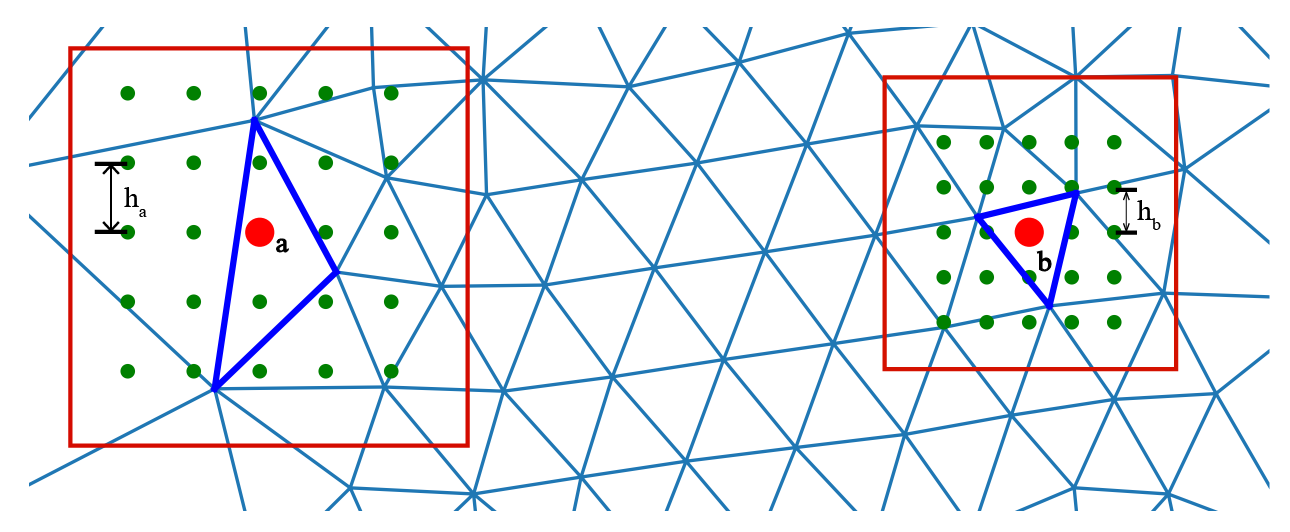}
    \caption{Example of $5 \times 5$ stencil for NNLCI local domain of dependence. Red dots: selected locations; green dots: local domain of dependence; blue lines: corresponding cells.}
    \label{fig:Stencils}
\end{figure}

First, we obtain the solutions from the simulations on the coarse and finer grids. Then, for each data location $(x,y)$, the local domain of dependence is determined. The intention is to filter the data at a proper scale to include all the local features for accurate prediction, while discarding far end information for low training cost. Fig.~\ref{fig:Stencils} gives two examples of this process. For the selected locations (red dots), and a given grid, the corresponding cells E are located (blue lines) can be found efficiently by adopting a hierarchical data structure described in Ref.~\cite{25}. The computational domain is divided into a binary tree of blocks of cells based on the cell centroid locations. By comparing the data location $(x,y)$ with the cell centroid, one can descend from the root to sub-blocks and find the cell $E$ containing it. This greatly reduces the computational cost and search time. The adjacent cells $E_{adj}$ that share edges or vertices with $E$ are then identified based on the connectivity information. Next, the local cell sizes are calculated. Here, the local cell size $h_E$ of a triangular cell $E$ is defined as the average cell size of itself and its adjacent cells:
\begin{equation}
    h_E =  \frac{1}{N} \sum_{k \in E_{adj}} \sqrt{A_k}
\end{equation}
where $A_k$ is the area of cell $k$, $E_{adj}$ is the adjacent cells to cell $E$, and $N$ is the total number of $A_k$ in the summation. We use a $5 \times 5$ rectangular stencil in 
$$[x-2h_{coarse},x+2h_{coarse}]\times[y-2h_{coarse},y+2h_{coarse}]$$
based on the coarse grid (where $h_{coarse}$ is the local cell size of the coarse grid near $(x,y)$) as the local domain of dependence, shown by green points. It should be noted that the size of the local domain of dependence varies with the local cell size. Also, during the process, some of the data points near the boundaries are discarded, since the selected stencil points are out of bound.

Once a stencil is determined, the values of state variables at the $5 \times 5$ locations determined through interpolation from the coarse and finer meshes are used as the input of a neural network. Note that the input also includes $h_{coarse}$ and $h_{finer}$, the local cell size at $(x,y)$ for the coarse and finer mesh. These two values are critical for the neural network to approximate the non-uniform properties, given the states interpolated from the coarse and finer meshes. The output is the predicted state variables at $(x,y)$. The procedure for determining the values of state variables is as follows. Let $\bm{\vec{x}}_s=(x_s,y_s)$ be one of the $5 \times 5$ stencil points in a particular stencil. We first locate the corresponding cell $E_0$ that contains this point using the hierarchical approach as above. The cell center of cell $E_0$ is denoted as $\bm{\Vec{x}}_0=(x_0,y_0)$ and the state variables of the cell is $\textbf{u}_0$. Then, the three adjacent triangles $E_1$, $E_2$ and $E_3$ are located, with cell centers and state as $(x_i,y_i)$ and $\textbf{u}_i$, where $i=1,2$ and $3$. The state at the stencil point location $(x_s,y_s)$ can be interpolated with a linear polynomial:
\begin{equation}
    \bm{u}(\bm{\vec{x}}_s-\bm{\Vec{x}}_0) = \bm{a}_0+\bm{a}_1(x_s-x_0)+\bm{a}_2(y_s-y_0)
\end{equation}
It is easily seen that $\bm{a}_0=\bm{u}_0$. To solve the coefficients $\bm{a}_1$ and $\bm{a}_2$, three combinations of cells can be used: $(E_0,E_1,E_2)$, $(E_0,E_2,E_3)$ and $(E_0,E_1,E_3)$. Each will determine a set of candidate values for $\bm{a}_1$ and $\bm{a}_2$, say ${\bm{a}_1}_k$ and ${\bm{a}_2}_k$. Then, the minmod function is applied to determine the best coefficient values:
\begin{equation}
    \bm{a}_i = minmod({\bm{a}_1}_k,{\bm{a}_2}_k),   k=1,2,3
\end{equation}
where the minmod function is defined as:
\begin{equation}
    minmod(a_1,a_2,...,a_n) = 
    \begin{cases}
        min(a_1,a_2,...,a_n), &\text{if all } a_i > 0 \\
        max(a_1,a_2,...,a_n), &\text{if all }  a_i < 0 \\
        0, &\text{otherwise}
    \end{cases}
\end{equation}

The procedure is repeated on the coarse, finer, and high-fidelity solutions to generate the input and reference values for the neural network. The local cell sizes are used as the inputs of the neural network to include the local mesh size information, which gives us an input size of 202, including the four state variables at $5 \times 5$ locations and the two local cell sizes. The output size is 4. The input-output pair is shown in Fig.~\ref{fig:Input-output}:
\begin{figure}[ht!]
    \centering
    \includegraphics[width=6in]{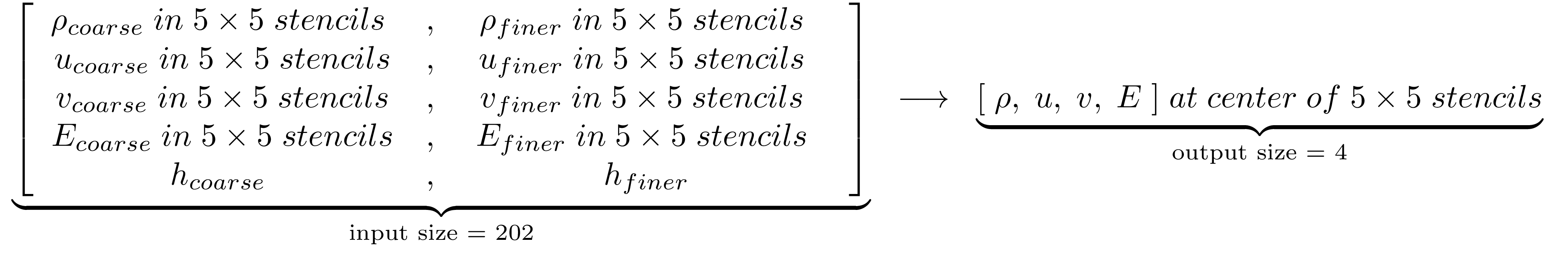}
    \caption{Input-output pair of the unstructured NNLCI. The input contains the state variable values on 5×5 stencil points, and the local mesh sizes of coarse and finer grids. The output is the state variable value at the center.}
    \label{fig:Input-output}
\end{figure}

\begin{table}[ht!]
\centering
\begin{tabular}{cc}
\hline
& Hyperparameters\\
\hline
Number of epochs & 50000\\
Number of hidden layers & 10\\
Network Structure & [202, $10 \times 500$, 4]\\
Learning rate & $1 \times 10^{-4}$ \\
$L_2$ regularization & $1 \times 10^{-8}$\\
Activation function & Tanh\\
Optimizer & Adam\\
\hline
\end{tabular}
\caption{\label{tab:nn_parameters} Hyperparameters of the neural network of local converging input (NNLCI)}
\end{table}
Table~\ref{tab:nn_parameters} lists the hyperparameters of the selected neural network, determined from manual search, in the present study. A network of 10 hidden layers of 500 nodes is designed to learn the mapping from low-fidelity solutions in a local domain to the high-fidelity solution at a point. The network is trained using the Adam optimizer with a learning rate of $1 \times 10^{-4}$ with $1 \times 10^{-8}$ regularization. The Tanh function is used as the activation function. The relative Mean squared error (RMSE) is selected as the loss function to measure the difference between NNLCI prediction and the true high-fidelity data:
\begin{equation}
     \mathcal{L} =  \frac{\sum_k \|\bm{\Tilde{u}_k}-\bm{u_k}\|_2^2}{\sum_k \|\bm{u_k}\|_2^2}
\end{equation}
where $\bm{\Tilde{u}_k}$ is the predicted results from NNLCI and $\bm{u_k}$ the data of high-fidelity simulation. For use on unstructured data, the RMSE needs to be weighted by the cell area. The cell-weighted RMSE is given by:
\begin{equation}
     \mathcal{L} = \frac{\sum_k  A_k\|\bm{\Tilde{u}_k}-\bm{u_k}\|_2^2}{\sum_k A_k\|\bm{u_k}\|_2^2}
\end{equation}

\section{Results and Discussion}
\label{Sec: Results and Discussion}
\begin{figure}[ht!]
    \centering
    \includegraphics[width=6in]{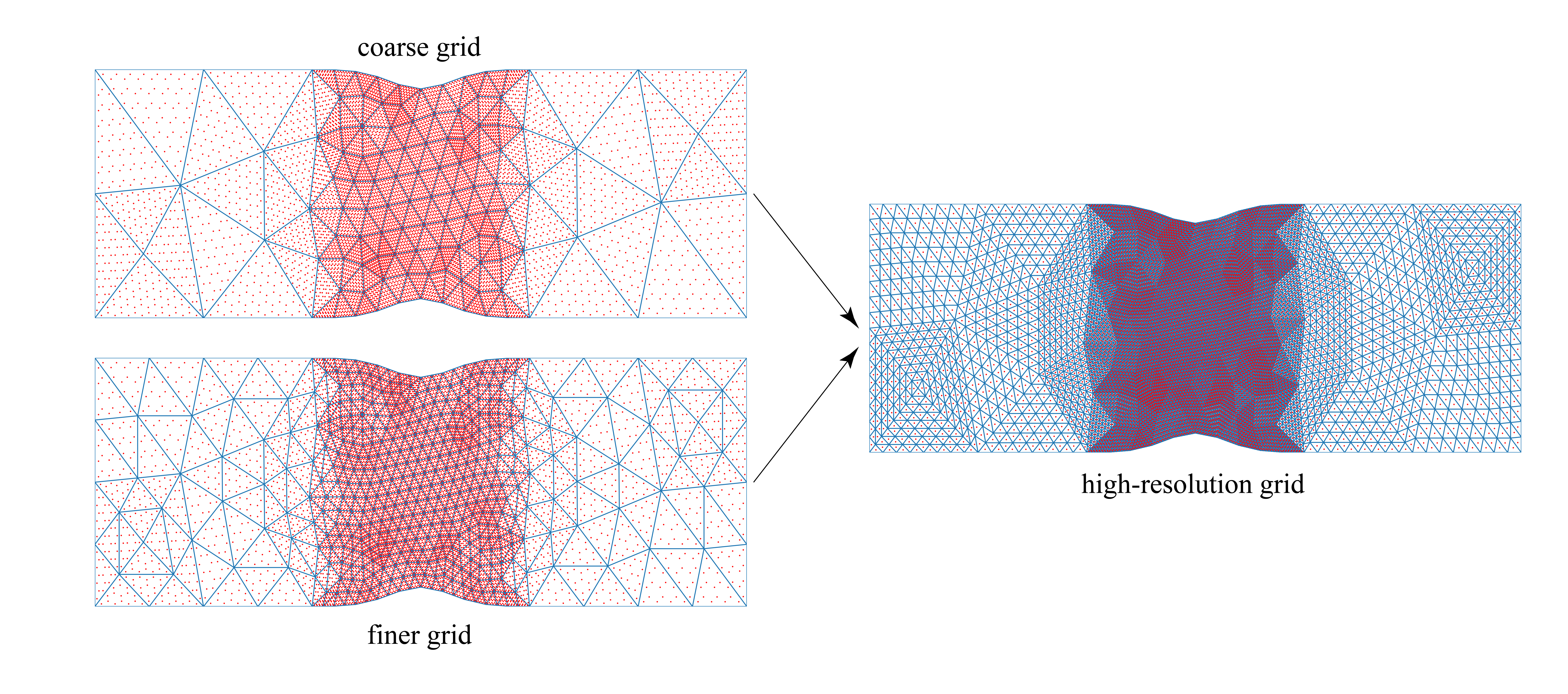}
    \caption{Example of the unstructured NNLCI prediction process. Data at prediction locations are interpolated on coarse and finer grids.}
    \label{fig:NNLCI_Pred_Inputs}
\end{figure}

IIn this section, we present the unstructured NNLCI predictions for several bump geometries. The coarse and finer solutions of each case are used as the inputs for the neural network. Unlike in the training cases, we select the cell centroids of the finest grids as the prediction locations to enhance the resolution of final prediction results. Fig.~\ref{fig:NNLCI_Pred_Inputs} shows an example of the input-image pairs. The red dots denote the selected locations for prediction. The interpolation technique described in Sec.~\ref{Sec: NNLCI} is used to get the data for these locations on the coarse and finer grids.

\begin{figure}[ht!]
    \centering
    \includegraphics[width=6in]{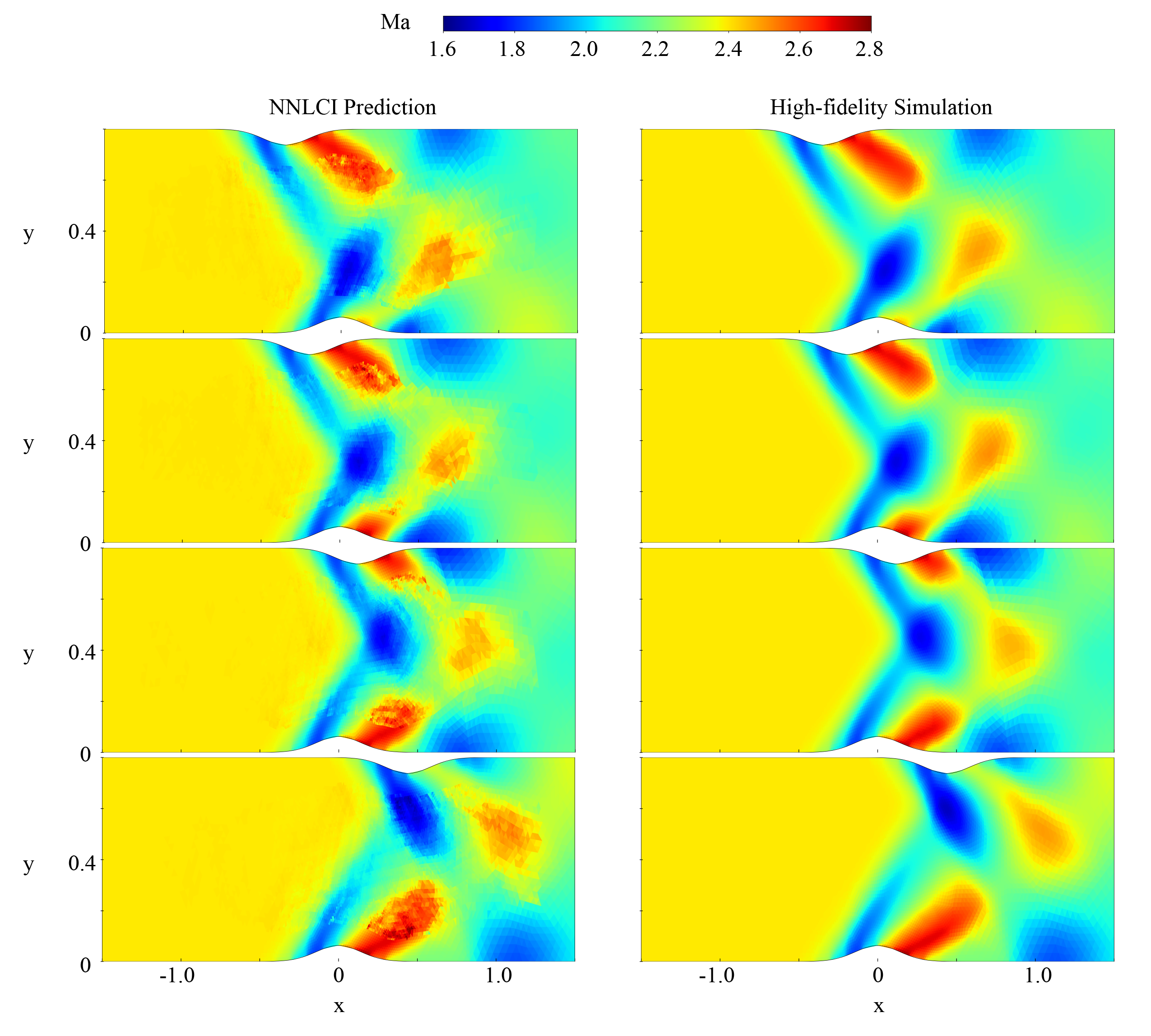}
    \caption{Mach number fields of the NNLCI prediction (left) and high-resolution simulation (right). Upper bump translations are $\Delta x=-0.35,-0.19,0.12,0.44$ respectively.}
    \label{fig:Mach_contour}
\end{figure}

Fig.~\ref{fig:Mach_contour} shows the NNLCI predicted Mach-number fields, for the upper bump translations of $\Delta x = -0.35$, $-0.19$, $0.12$ and $0.44$. High-fidelity simulation results are also presented for comparison. The flowfields exhibit different behaviors and features, depending on the bump translation. In the case of $\Delta x=-0.35$, the upper shock develops upstream and intersects with the lower shock near the lower bump. As a consequence, the flow expansion area is shifted toward the lower wall. The secondary shock is well-developed downstream of the upper bump, while it can hardly be observed near the lower bump. On the contrary, the case of $\Delta x=0.44$ has opposite behavior. For cases of $\Delta x=0.12$ and $-0.19$, a secondary shock is observed for both the upper and lower bumps, with shock intersection near the center of the channel. Regardless of the complex flow features, the NNLCI method achieves promising results. For all four cases, it accurately captures the primary shock location and structure. The complex nonlinear intersection of the upper and lower shocks is precisely predicted. The expansion region and secondary shock are well-reconstructed.

Table~\ref{tab:Mach error} summarizes the prediction error over the four cases. Here, the area-weighted relative $L_1$ norm is used as the measure for performance:
\begin{equation}
     \mathcal{L}_1 =  \frac{\sum_k  A_k\|\bm{\Tilde{u}_k}-\bm{u_k}\|_1}{\sum_k A_k\|\bm{u_k}\|_1}
\end{equation}
As comparison, the area-weighted relative root mean square error (RRMSE) is added:
\begin{equation}
     \mathcal{L} = \sqrt{ \frac{\sum_k  A_k\|\bm{\Tilde{u}_k}-\bm{u_k}\|_2^2}{\sum_k A_k\|\bm{u_k}\|_2^2}}
\end{equation}

\begin{table}[ht!]
\renewcommand{\arraystretch}{1.2} 
\centering
\begin{tabular}{cccc}
\hline
Bump Translation & Relative $L_1$ Norm & RRMSE & Low-fidelity Relative $L_1$ Norm\\\hline
-0.35 & $0.701\%$ & $1.122\%$ & $6.421\%$ \\
-0.19 & $0.738\%$ & $1.115\%$ & $6.345\%$ \\
0.12 & $0.757\%$ & $1.097\%$ & $6.813\%$ \\
0.44 & $0.449\%$ & $0.699\%$ & $6.647\%$ \\
Total & $0.659\%$ & $1.021\%$ & $6.554\%$\\
\hline
\end{tabular}
\caption{\label{tab:Mach error} Relative $L_1$ norm and relative root mean square error for prediction cases}
\end{table}

For all four prediction cases, the NNLCI method achieves a relative $L_1$ error around one percent. As a comparison, the relative $L_1$ error is measured on the low-fidelity simulation data. The data are interpolated using the same methods on the desired prediction locations and used as the inputs for the NNLCI. It can be seen that the low-fidelity error is around $6.5\%$, while the overall error of the NNLCI method is below $1\%$. The prediction accuracy is improved by about 10 times.

\begin{figure}[ht!]
    \centering
    \includegraphics[width=6in]{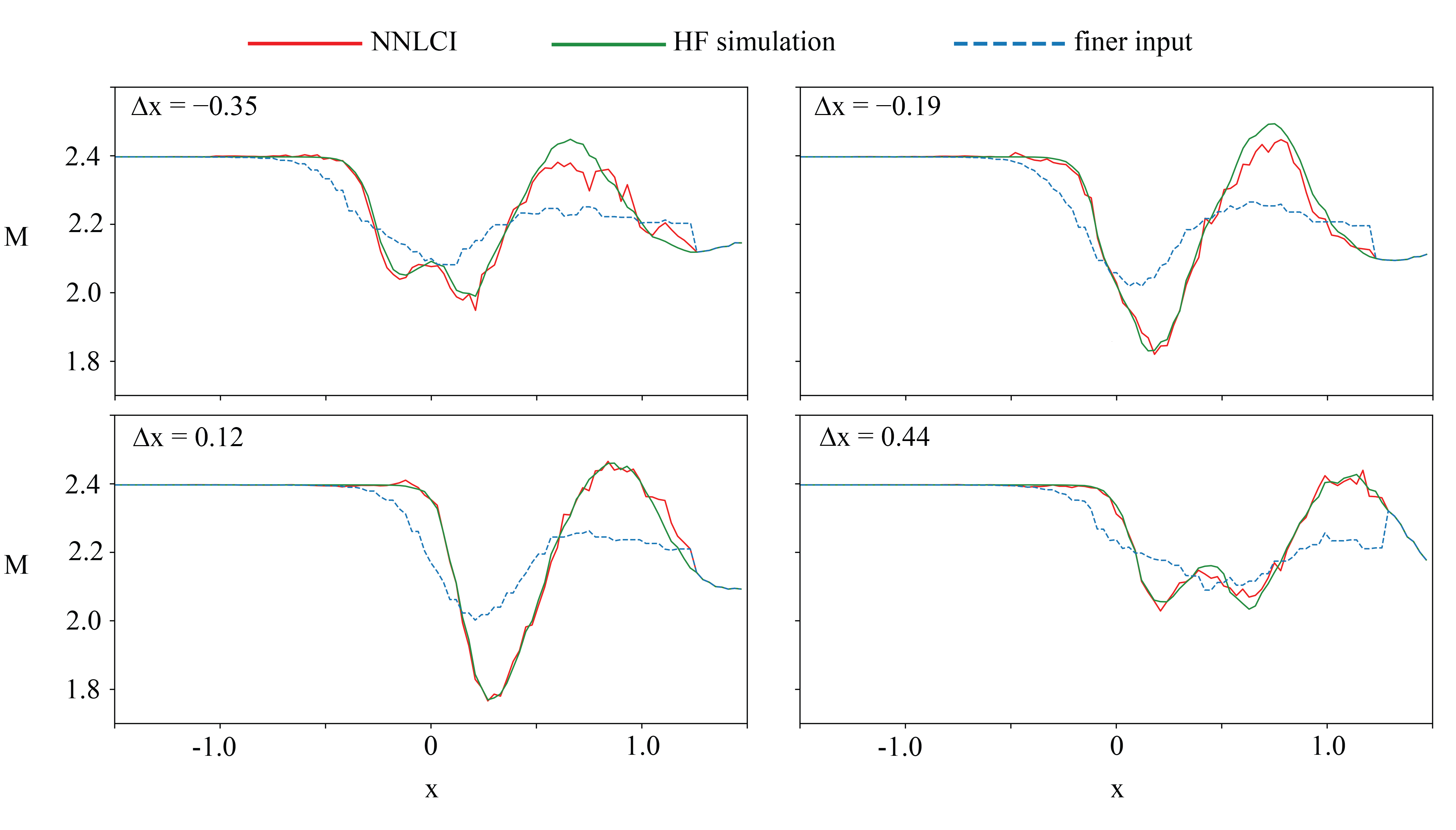}
    \caption{Mach-number distribution along the centerline of the channel $y=0.4$}
    \label{fig:Cross_section}
\end{figure}

Fig.~\ref{fig:Cross_section} shows the Mach-number distribution along the centerline of the channel, $y=0.4$. The NNLCI prediction (red line) agrees well with the high-fidelity simulation result (green line). The finer input from low-resolution simulation results is also presented (blue line). The flow development and shock interaction are accurately captured by the NNLCI prediction for all the test cases. The prediction of shock shape and location is greatly enhanced compared with the low-fidelity simulation. These results testify to the effectiveness of the NNLCI method on the flow prediction for the entire field, including regions with smooth evolution and shock discontinuities \cite{19,20}.

\begin{figure}[ht!]
    \centering
    \includegraphics[width=6in]{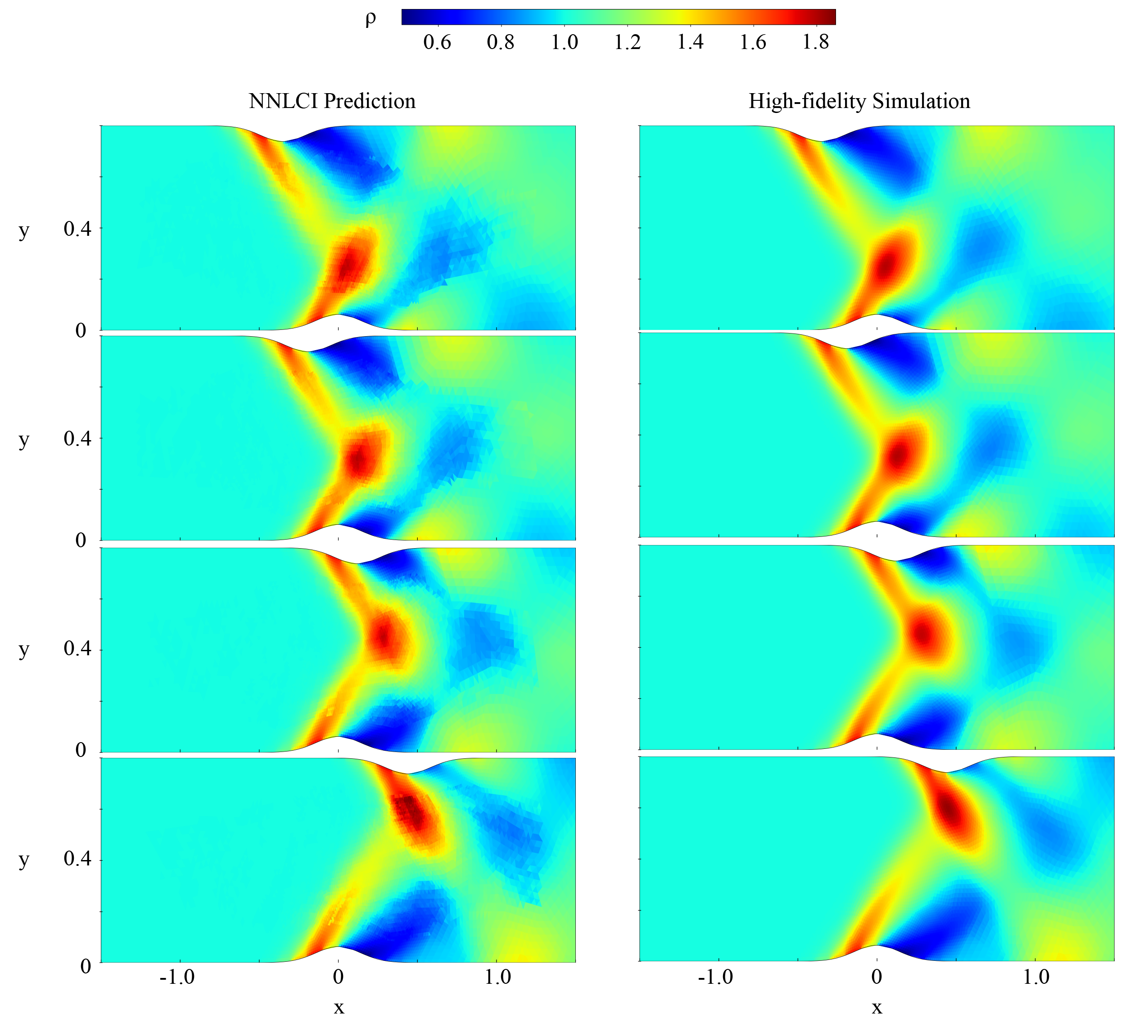}
    \caption{Density fields: NNLCI prediction (left) and high-fidelity simulation results (right). Upper bump translations are  $\Delta x=-0.35,-0.19,0.12,0.44$, respectively.}
    \label{fig:Density_contour}
\end{figure}

\begin{figure}[ht!]
    \centering
    \includegraphics[width=6in]{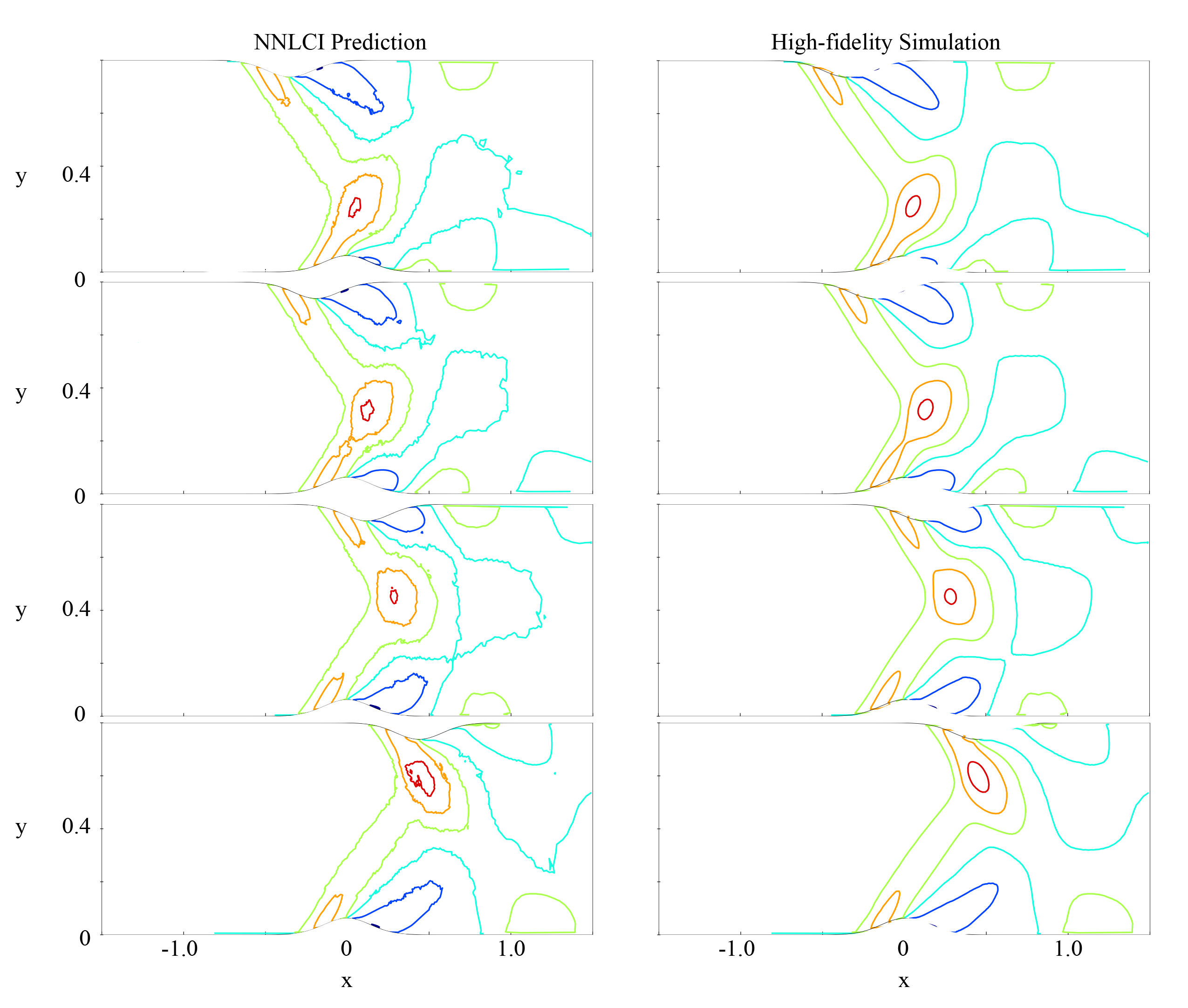}
    \caption{Density contour gradient: NNLCI prediction (left) and high-fidelity simulation results (right). Upper bump translations are  $\Delta x=-0.35,-0.19,0.12,0.44$, respectively.}
    \label{fig:Density_gradient}
\end{figure}

\begin{figure}[ht!]
    \centering
    \includegraphics[width=6in]{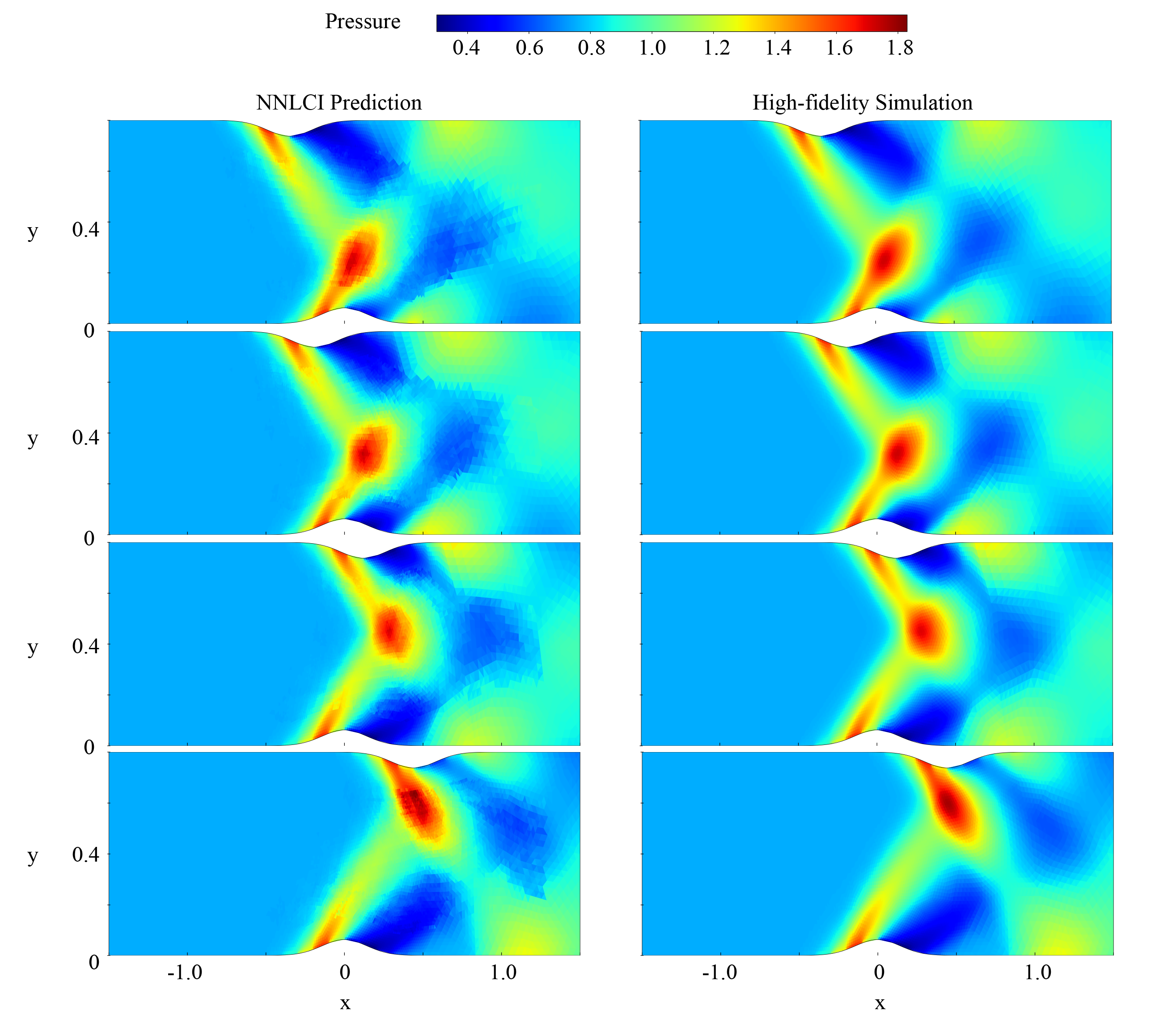}
    \caption{Pressure fields: NNLCI prediction (left) and high-fidelity simulation results (right). Upper bump translations are  $\Delta x=-0.35,-0.19,0.12,0.44$ respectively.}
    \label{fig:Pressure_contour}
\end{figure}

As introduced in Sec 2.3, the NNLCI network can predict the four state variables simultaneously. It is thus interesting to examine the NNLCI performance on each state variable separately. Figs.~\ref{fig:Density_contour}, ~\ref{fig:Density_gradient} and \ref{fig:Pressure_contour} show the results of the density fields, density contour gradient and pressure fields. The NNLCI method can accurately capture the shape and magnitude of discontinuities for both variables. In addition, the prediction of smooth regions closely matches the simulation results. The NNLCI method achieves an accuracy of more than $99\%$ for all cases. For a new bump configuration, the NNLCI method can predict all the state variables with one neural network, eliminating the need for repeated training of multiple networks for different variables. 

Compared with the physics-informed methods, the NNLCI method can capture shock and shock interactions accurately. In addition, different from the conventional global-to-global deep learning methods, the NNLCI method avoids the use of complex neural network structure and eliminates the need for large dataset with small training gap. This greatly reduces the training time and computational cost. The prediction accuracy of local patches is improved with richer details. For the present study, the wall time for low and high-fidelity simulation of each design setting takes 10 seconds and 30 minutes respectively, on a single CPU (Intel Core i7-10750H). In comparison, the NNLCI is able to predict a new case in less than 1 second on the same hardware. The time saving is more than 100 times. 

\begin{figure}[!ht]
    \centering
    \begin{subfigure}{4in}
    \centering
        \includegraphics[width=4in]{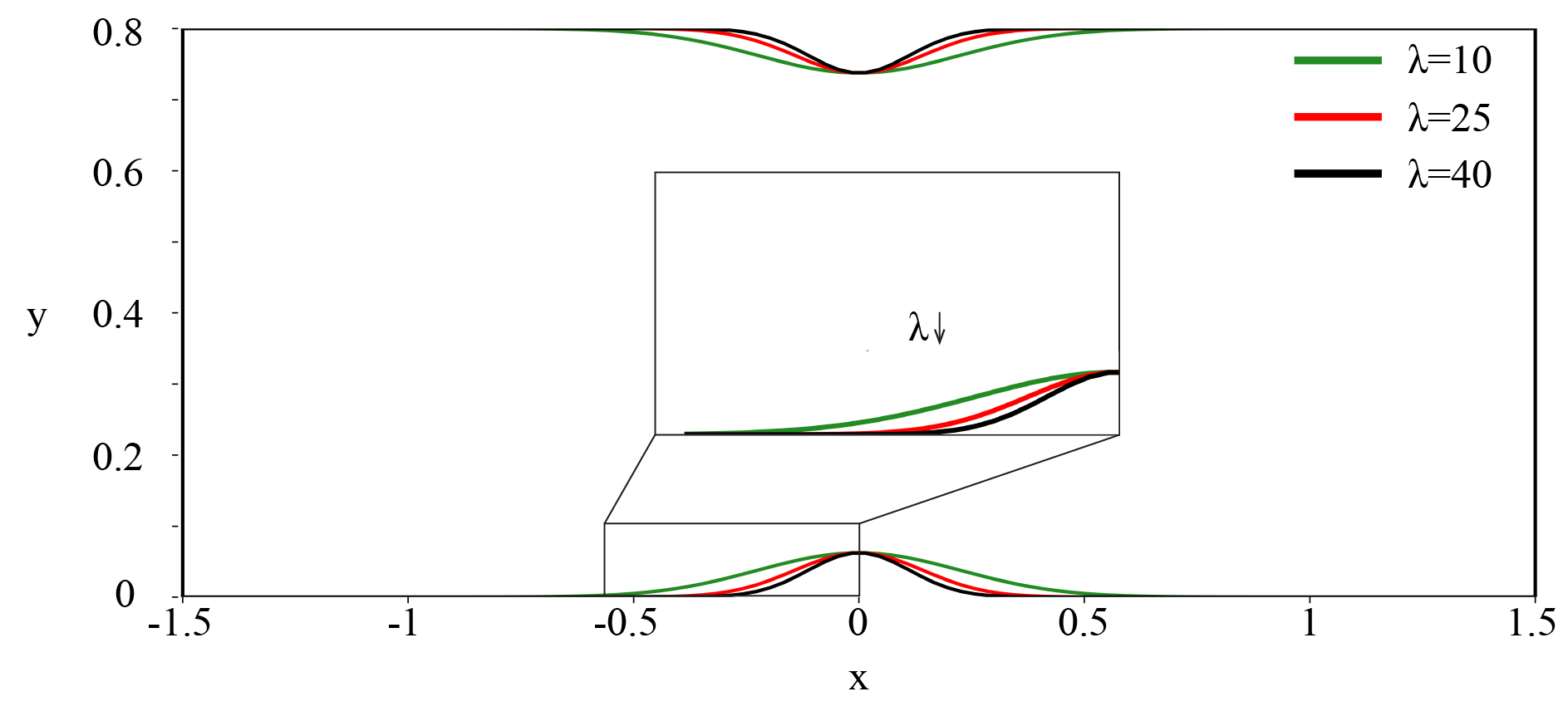}
        \label{t=0}
    \end{subfigure}
    
    \begin{subfigure}{4in}  
    \centering
        \includegraphics[width=4in]{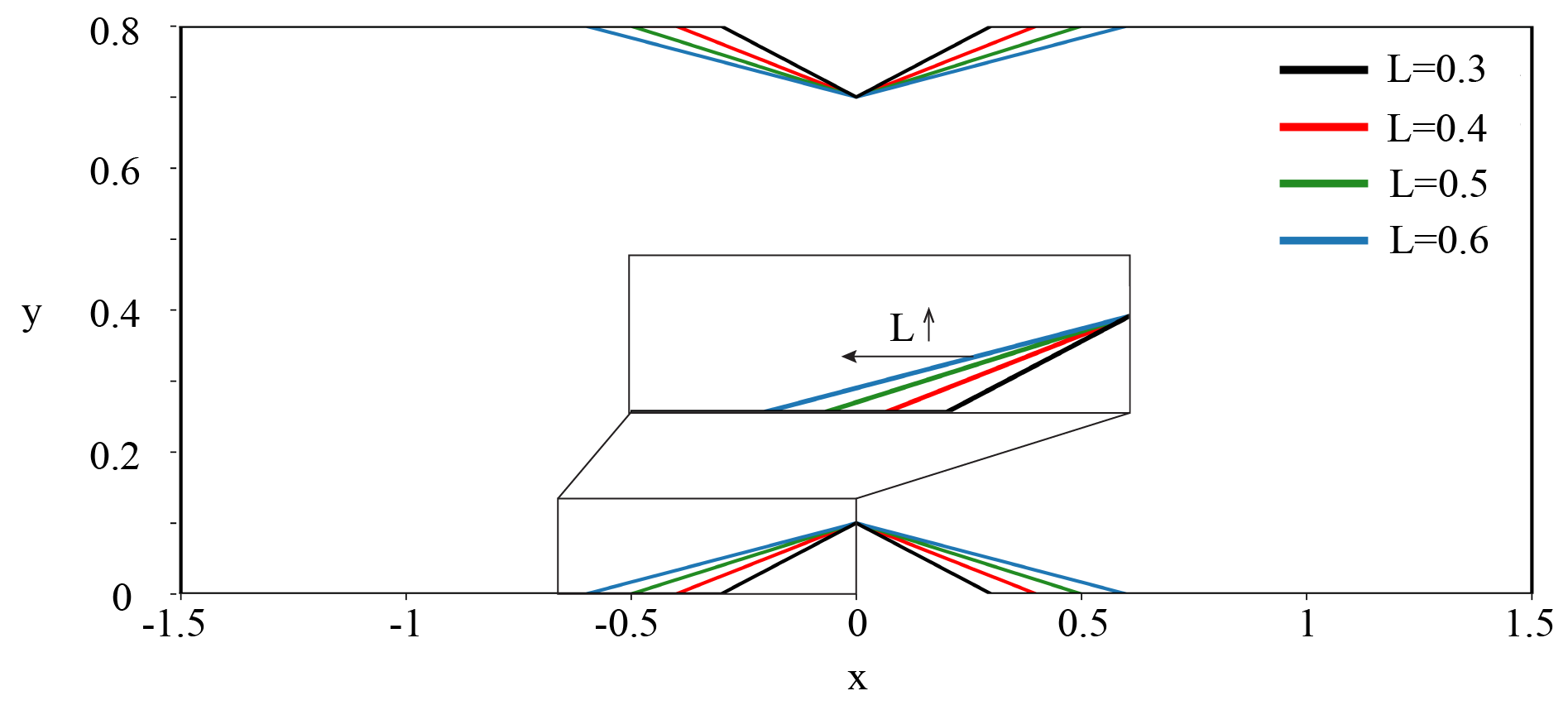}
        \label{t=180}
    \end{subfigure}
    \caption{Gaussian and triangular bumps with parameters $\lambda$ and $L$.}
    \label{fig:bump_shape}
\end{figure}

To further evaluate the effectiveness and flexibility of the NNLCI method, two different bump geometries are investigated: a triangular bump and a Gaussian bump with changing variance, as shown in Fig.~\ref{fig:bump_shape}. The shape of the Gaussian bump is varied by tuning the variance $\lambda$ of the governing function, with the amplitude fixed.
\begin{equation}
    y = 0.0625e^{-\lambda x^2}.
\end{equation}
Second, the Gaussian bump is replaced by the triangular wedge. The height of the wedge $h$ is fixed as $0.1$, while the length of the wedge $L$ is changed in the range of $[0.3, 0.6]$ with increment of $0.1$. For both cases, the angle of attack varies with the bump shape, and the flow exhibits different behaviors.

\begin{figure}[ht!]
    \centering
    \includegraphics[width=6in]{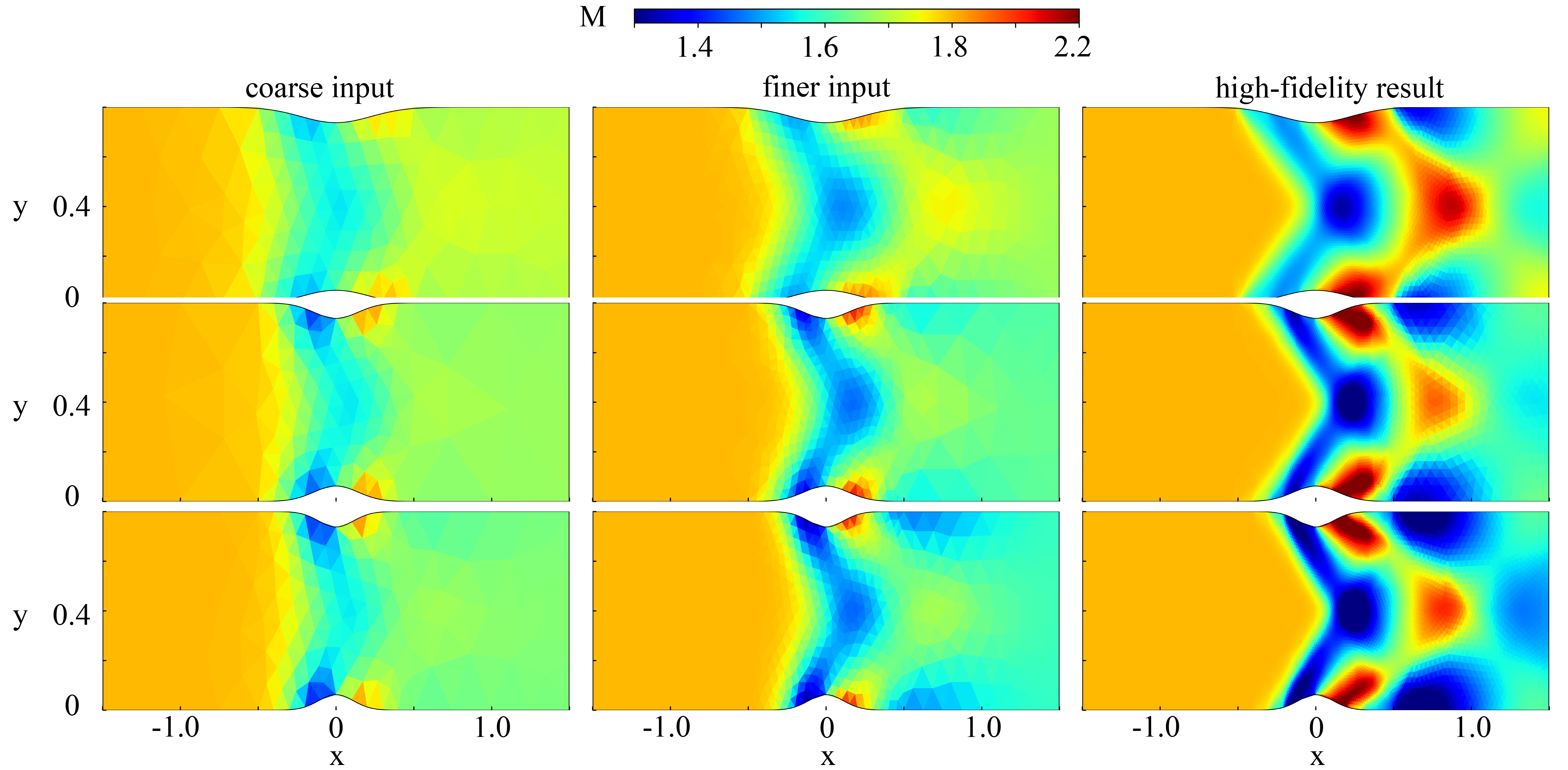}
    \caption{Simulation results of Mach-number fields for Gaussian bumps with $\lambda = 10$, $25$, and $40$. The coarse, finer inputs, and high-fidelity results are shown, from left to right.}
    \label{fig:G_Simulation}
\end{figure}

\begin{figure}[ht!]
    \centering
    \includegraphics[width=6in]{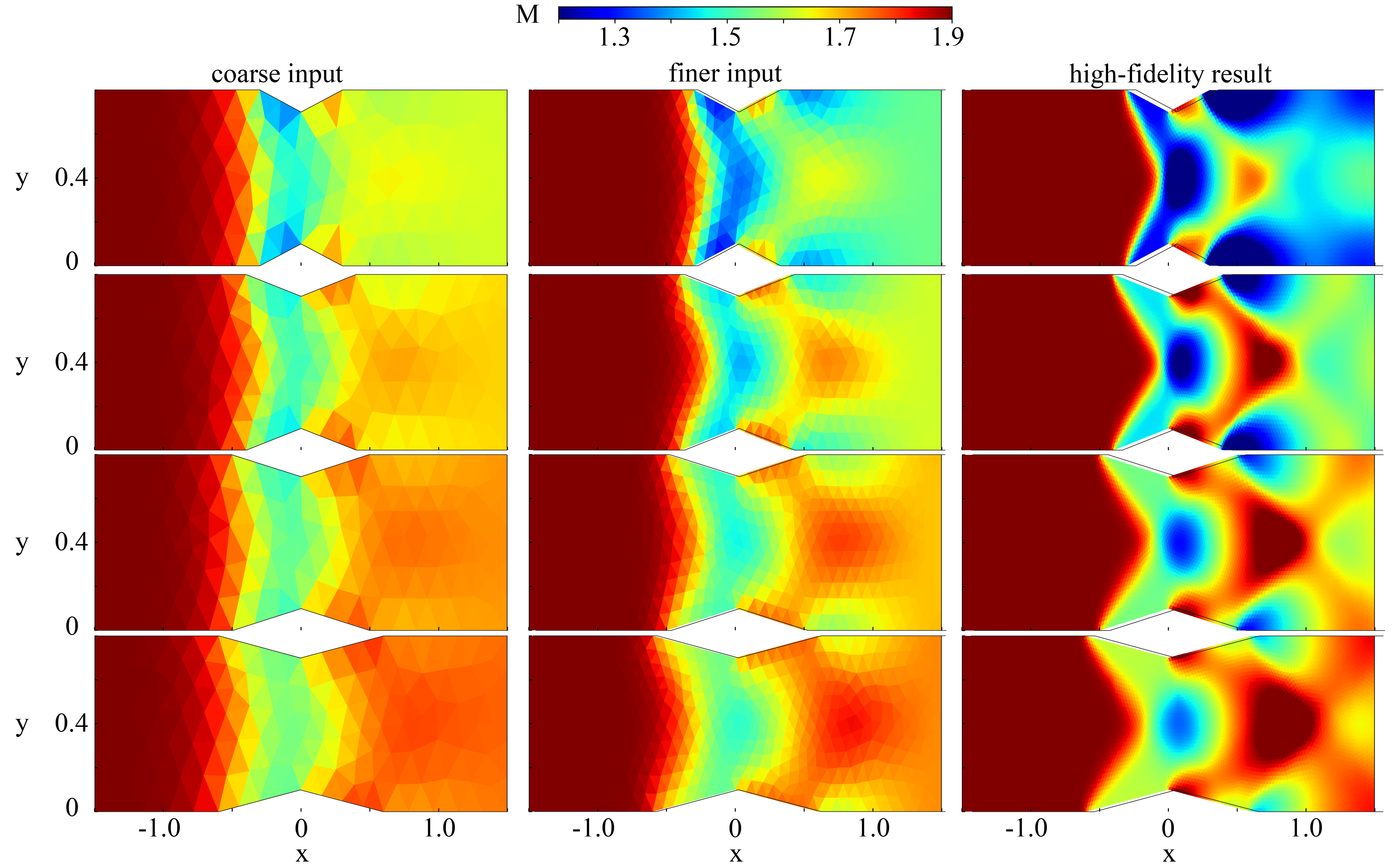}
    \caption{Simulation results of Mach-number fields for triangular bumps with $L = 0.3$, $0.4$, $0.5$ and $0.6$. The coarse, finer inputs, and high-fidelity results are shown, from left to right.}
    \label{fig:Tri_Simulation}
\end{figure}

Figs.~\ref{fig:G_Simulation} and \ref{fig:Tri_Simulation} show numerical simulation results of the Mach number fields for the two different kinds of bump shapes. Tuning the geometric parameter $\lambda$ or $L$ changes the angle of attack of the incoming flow, and the shock structure and dynamics vary as a result. With the triangular bump, as the angle of attack decreases, the secondary shock structure moves downstream and the Mach number alters significantly. Such enormous change of flow behavior greatly increases the difficulty of prediction. In the present study, instead of constructing a separate neural network, the new cases are trained together with the bump translation dataset. The resultant NNLCI is expected to predict the flowfield for all bump geometries with high accuracy simultaneously. Table~\ref{tab:shape_table} summarizes the training and validation cases. For each bump shape, the free-stream Mach number $M_{infty}$ is perturbed by $\pm 5\%$

\begin{table}[ht!]
\renewcommand{\arraystretch}{1.5} 
\centering
\begin{tabular}{cccc}
\hline
Set & Bump Type & Description & $M_{\infty}$ Perturbation\\
\hline
\multirow{2}{*}{Training} & Gaussian Bump & $\lambda = 10, 25, 40$ & \multirow{2}{*}{$M_{\infty} = 0, \pm 5\%$} \\
                     & Triangular Wedge & $L=0.3,0.4,0.5,0.6$ & \\
\hline
\multirow{2}{*}{Testing} & Gaussian Bump & $\lambda = 28$ & \multirow{2}{*}{None}\\
                     & Triangular Wedge & $L=0.38$ & \\
\hline
\end{tabular}
\caption{\label{tab:shape_table}Training and validation cases for Gaussian and triangular bumps.}
\end{table}

\begin{figure}[ht!]
    \centering
    \includegraphics[width=6in]{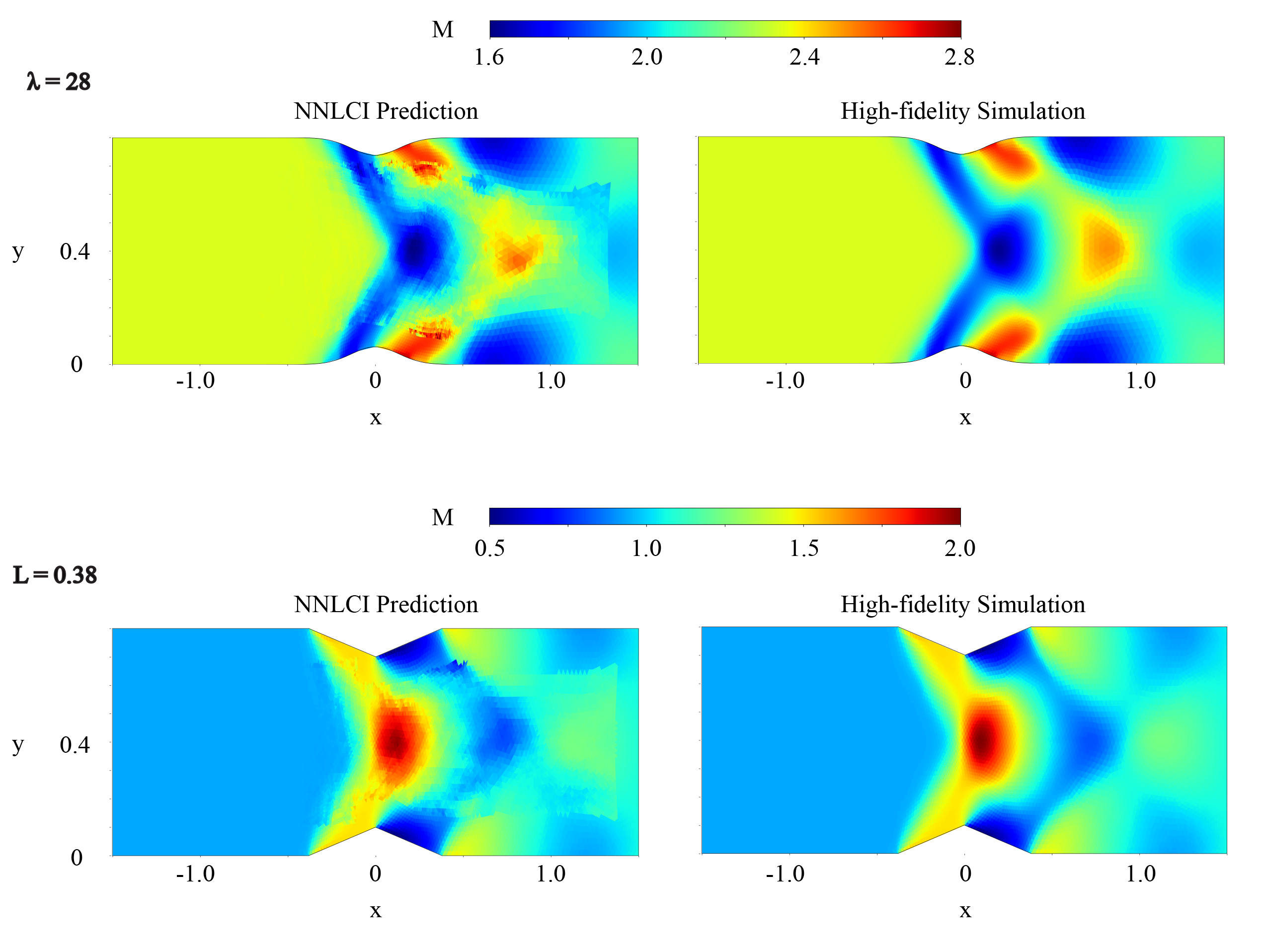}
    \caption{Mach-number fields calculated by the NNLCI method (left) and high-fidelity simulation (right) for Gaussian bump $\lambda = 28$ (upper) and triangular bump $L=0.38$ (lower).}
    \label{fig:Mach_contour_shape}
\end{figure}

Fig.~\ref{fig:Mach_contour_shape} shows the Mach-number fields for the NNLCI predictions and high-fidelity simulation. As seen in the plot, the shock develops from different locations with variation of the bump length. For example, the wedge bump case has a larger oblique shock angle, due to the large angle of attack. As a consequence, the subsonic region is compressed and the secondary shock strength is weaker, due to insufficient expansion. A well-developed smooth region is observed further downstream. The Gaussian bump result exhibits opposite behavior, but the NNLCI method accurately predicts the shock structure and the smooth region behavior for both cases. The prediction error is summarized in table~\ref{tab:Mach error_shape}. In both cases, the NNLCI achieves an accuracy of more than $98 \%$, while maintaining precision for the bump translation cases. Compared to the low-fidelity simulation results, the NNLCI is able to increase precision by three times.

As noted above, only one NNLCI is built and trained for all the cases. This avoids the construction and training of new neural networks for different geometries. In practice, the NNLCI method can greatly reduce design turnaround time; for the prediction of a new design feature, one can train the existing NNLCI model with supplementary data and obtain results in a short time. Overall, this further validates the effectiveness of NNLCI method for real applications.

\begin{table}[ht!]
\renewcommand{\arraystretch}{1.2} 
\centering
\begin{tabular}{cccc}
\hline
Bump Shape & Relative $L_1$ Norm & RRMSE & Low-fidelity Relative $L_1$ Norm\\\hline
$\lambda=28$ & $1.279\%$ & $1.978\%$ & $6.389\%$ \\
$L=0.38$ & $2.055\%$ & $3.123\%$ & $7.125\%$ \\
Bump translation & $1.141\%$ & $1.783\%$ & $6.554\%$ \\
Total & $1.352\%$ & $2.176\%$ & $6.624\%$\\
\hline
\end{tabular}
\caption{\label{tab:Mach error_shape} Relative $L_1$ norm and relative root mean square error for various bump shapes.}
\end{table}

\begin{figure}[ht!]
    \centering
    \includegraphics[width=6in]{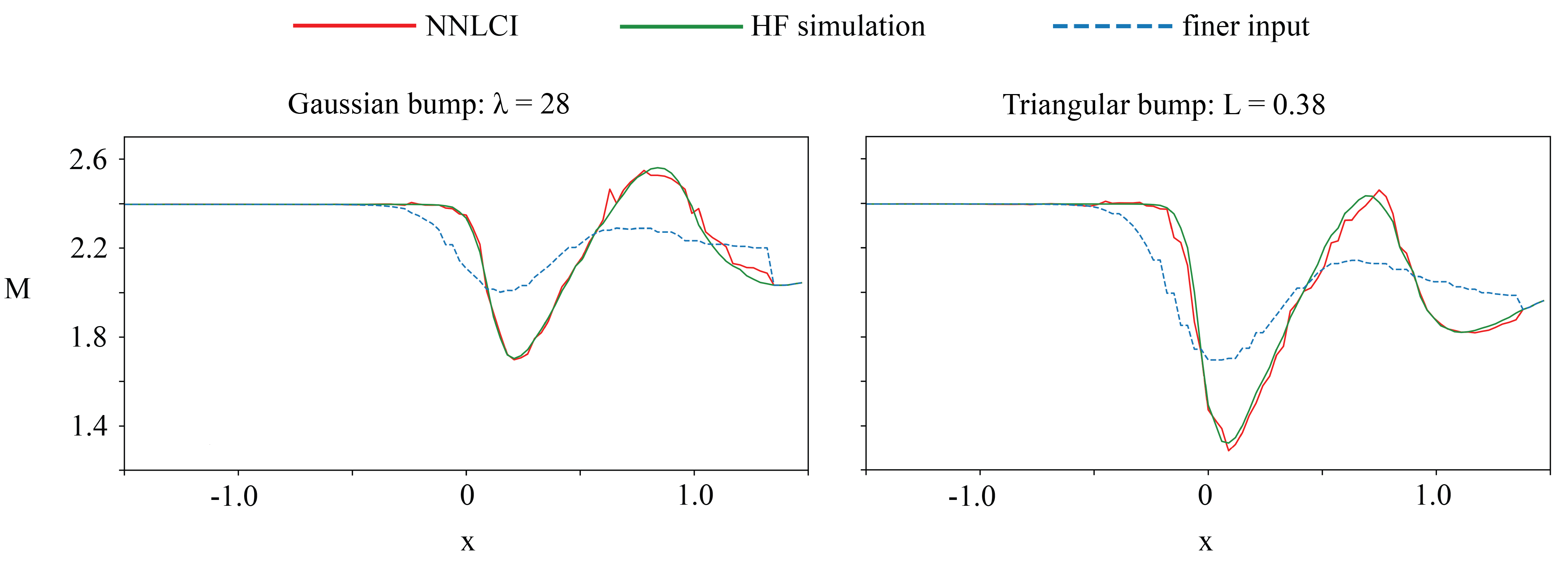}
    \caption{Cross-sectional plot of the Mach number contour at the center of the channel $y=0.4$ for Gaussian bump case $\lambda = 28$ (left) and triangular bump $L=0.38$ (right).}
    \label{fig:Shape_Cross_section}
\end{figure}

\begin{figure}[ht!]
    \centering
    \includegraphics[width=6in]{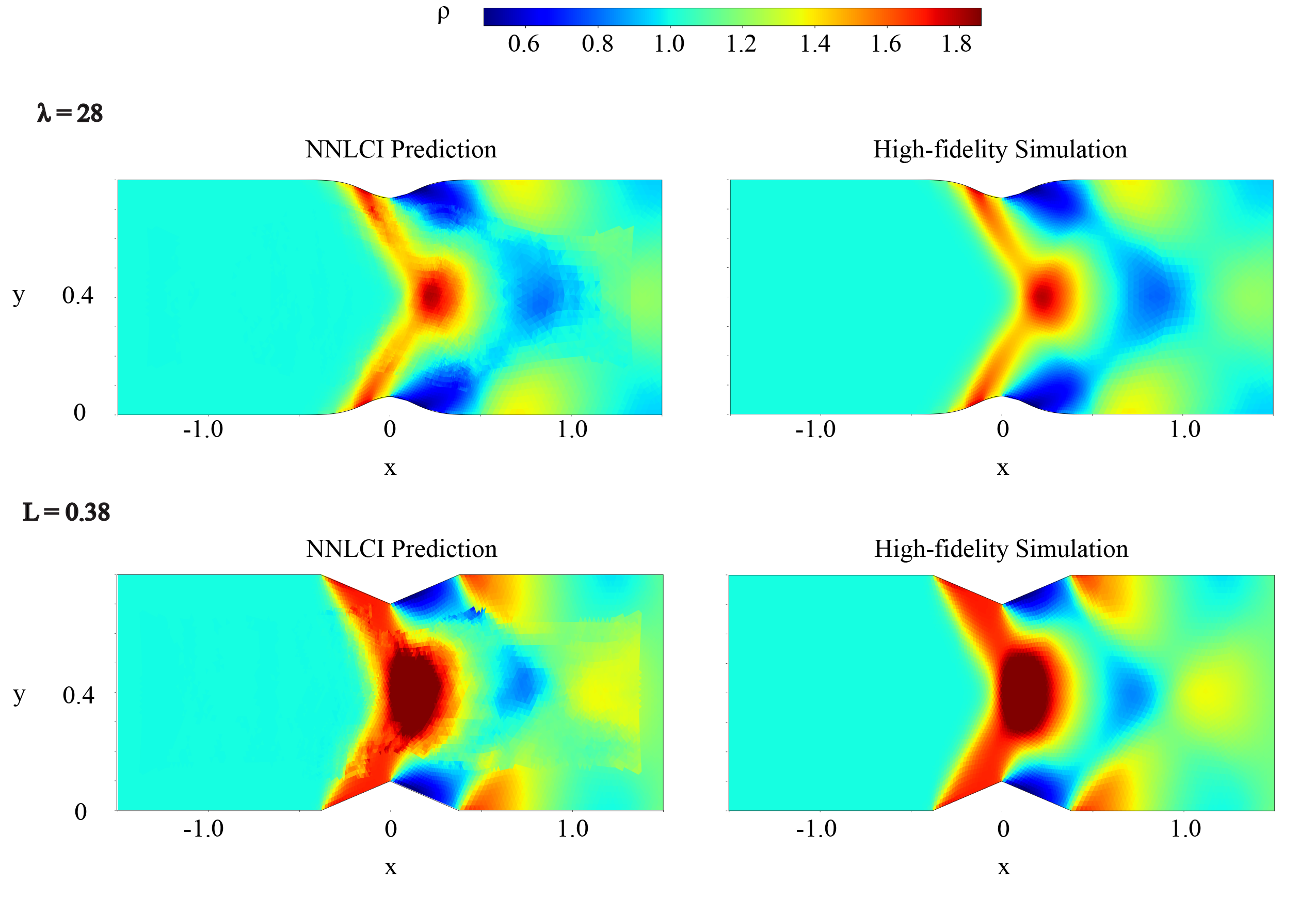}
    \caption{Density fields: NNLCI prediction (left) and high-fidelity simulation results (right) for Gaussian bump case $\lambda = 28$ (upper) and triangular bump $L=0.38$ (lower).}
    \label{fig:Shape_Density}
\end{figure}

\begin{figure}[ht!]
    \centering
    \includegraphics[width=6in]{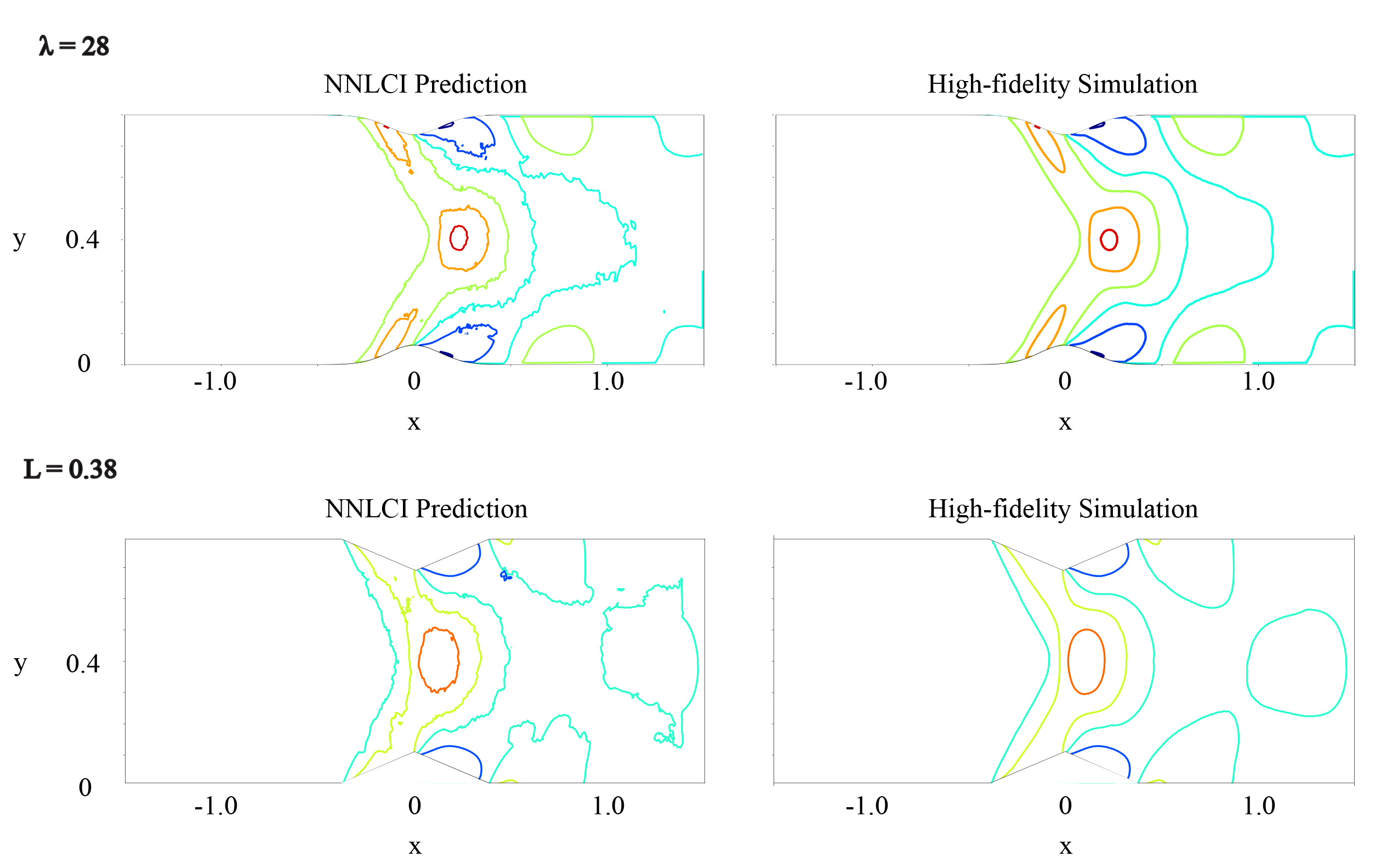}
    \caption{Density contour gradient: NNLCI prediction (left) and high-fidelity simulation results (right) for Gaussian bump case $\lambda = 28$ (upper) and triangular bump $L=0.38$ (lower).}
    \label{fig:Shape_Density_gradient}
\end{figure}

\begin{figure}[ht!]
    \centering
    \includegraphics[width=6in]{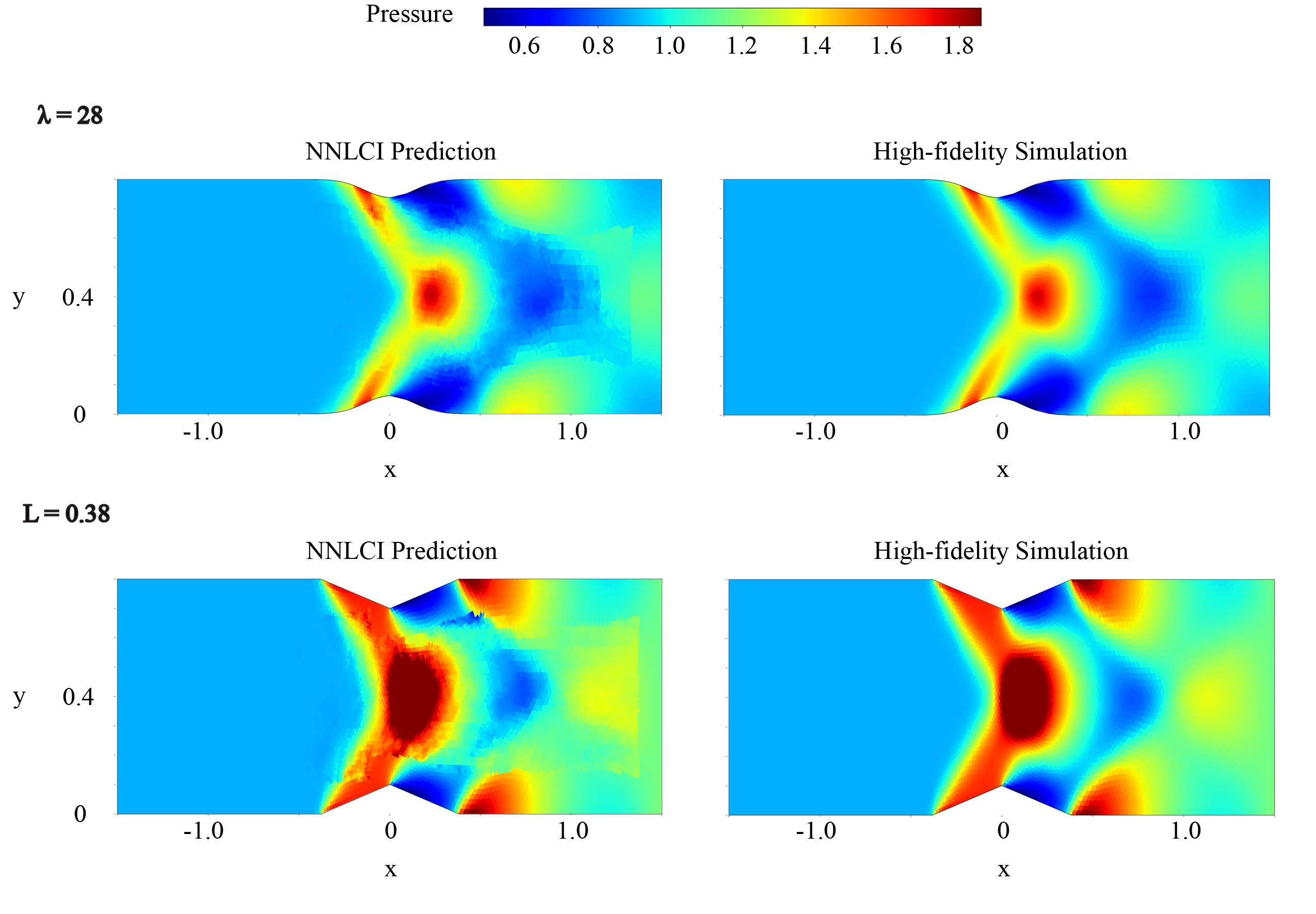}
    \caption{Pressure fields: NNLCI prediction (left) and high-fidelity simulation results (right) for Gaussian bump case $\lambda = 28$ (upper) and triangular bump $L=0.38$ (lower).}
    \label{fig:Shape_Pressure}
\end{figure}

A cross-sectional plot of Mach number contour is shown in Fig.~\ref{fig:Shape_Cross_section}. For the two cases with different bump shapes, the NNLCI is able to accurately capture the shock location and amplitude, despite the large deviation in flow behavior. In addition, Figs.~\ref{fig:Shape_Density}, ~\ref{fig:Shape_Density_gradient} and \ref{fig:Shape_Pressure} show density fields, density contour gradient and pressure fields for the two prediction cases. Despite the deviation near the boundary of the upper bump due to lack of data, the NNLCI can capture the development and variation of the shock behavior to high accuracy.

\section{Conclusion}
\label{Sec: conclusion}
This paper presents a neural network with local converging inputs (NNLCI) for unstructured data. The proposed model employs a novel sampling and interpolation technique to construct local converging inputs from low-resolution simulation results. The neural network builds up the regression from the local converging inputs to the high-fidelity prediction results.
To demonstrate, the unstructured NNLCI is applied to predict supersonic inviscid flow in a channel with two bumps. The upper bump geometry is perturbed to create complex shock intersection structures. A detailed comparison between the unstructured NNLCI prediction and the high-resolution simulation results is carried out. The NNLCI is able to accurately capture the shock structure. The density and pressure profiles are examined to demonstrate the capabilities of NNLCI in predicting multiple flow variables simultaneously. Furthermore, the unstructured NNLCI is successfully applied to other bump geometries. Without the training of a new model, the NNLCI can absorb additional features from new training data and produce accurate predictions of new design geometries.

\clearpage
\bibliography{mybibfile}

\end{document}